\documentclass[english]{article}
\usepackage[T1]{fontenc}
\usepackage[latin9]{inputenc}
\usepackage{amsmath}
\usepackage{graphicx}
\usepackage{amssymb}
\usepackage{esint}



\usepackage{babel}

\begin{document}

\title{Introduction to generalised C\'{e}saro convergence III}

\author{Richard Stone}
\maketitle
\begin{abstract}
\begin{eqnarray*}
\begin{array}{cc}
Beauty \; is \; its \; own \; excuse \; for \; being \\
\end{array}
\end{eqnarray*}

This is the third and last of three papers introducing generalised C\'{e}saro convergence and is split into two parts. In part 1 we introduce the notion of a "C\'{e}saro-adapted scale" and use it to prove the key generalised C\'{e}saro summation/convergence theorems developed in [1]. We also use it to trivially extend these results to the case of remainder C\'{e}saro summation/convergence relative to arbitrary $z_{0}\in\mathbb{C}$ (not just $z_{0}=0$). In the course of the working we introduce the concepts of "formal symbols" and "formal function elements", which allow us to express many results in extremely compact form and simplify our arguments considerably. 

Part 2 is self-contained and devoted to further exploring this "formal" world. We express a number of additional results in surprisingly compact form using formal symbols and function elements, and use them to give simple proofs of several non-trivial results. We also investigate their fascinating properties. These include the need to avoid evaluating too early; the consequent need to retain stand-alone zeros (both "to the left" and "to the right") lest they be brought back to life before evaluation; and the need to use continuous limits to resolve singular ratios in final evaluation when required. Finally we consider in detail the formal extension we have introduced of our C\'{e}saro-adapted scale to a 1-parameter continuum of period-1 functions $\overset{\lor}{q}_{\rho}(\alpha)$, $\rho\in\mathbb{C}$. We analyse their distributional aspects when $\rho\in\mathbb{Z}_{<0}$ and derive their Fourier-series coefficients in general. We conclude with a miscellany of further observations, including a formal re-casting of the general Euler-McLaurin sum formula in very compact form, and a number of additional analytical and combinatorial characteristics of the $\overset{\lor}{q}_{\rho}(\alpha)$ and associated operators.
\end{abstract}

\section{Introduction to Part 1}

Part 1 of this paper is brief and fulfils promises made in [1] by proving and then extending two key results developed there regarding the concept of generalised C\'{e}saro convergence.

In section 2.1 we first prove theorem 4 from [1] regarding the generalised C\'{e}saro limits of expressions of the form $X^{\delta}\alpha^{r}$ and $k^{\delta}\alpha^{r}$. We derive these as \textit{strong} C\'{e}saro asymptotic results by introducing the notion of a "C\'{e}saro-adapted scale" - a collection of related functions with properties well-adapted to iterative application of the C\'{e}saro operator, $P$.

In this case, the functions in this scale are closely-related to the period-1 versions of the famous Bernoulli polynomials, and we cover this relationship in section 2.2, where we also define various alternative normalisations of them. Just as the Bernoulli numbers can be defined via powers of a formal symbol, $B$, we also show in section 2.2 that our preferred normalised set, $\overset{\lor}{q}_{n}(\alpha)$, have a remarkably compact expression in terms of a related formal symbol, $\tau$, which we introduce. In turn this expression allow us generalise our definition of the $\overset{\lor}{q}_{n}(\alpha)$ to a full 1-parameter family of functions $\overset{\lor}{q}_{\rho}(\alpha)$ for arbitrary $\rho\in\mathbb{C}$ in a natural way. Finally in section 2.2 we further introduce the concept of a formal function element, since using this formalism will enable us to give a beautiful, compact C\'{e}saro expression for the remainder term in the p-sum $\sum_{j=1}^{k}j^{-s}$, and this will be critical in proving the next result in this paper.

This next result is theorem 3 from [1], giving the strong C\'{e}saro asymptotic relationship for $\sum_{j=1}^{k}j^{-s}$. This is perhaps the main result of that paper and it underlay the detailed analysis of the Riemann zeta function and its singularities within a C\'{e}saro framework undertaken there. In proving this result in section 2.3 we first derive the compact expression just mentioned and then proceed by invoking theorem 4 just proven. As part of the proof we are able immediately to derive the explicit form of the regular polynomial in $P$ involved in obtaining the C\'{e}saro limit of the p-sum at any given $s\in\mathbb{C\setminus}\{1\}$. We conclude section 2.3 with two brief comments, one of which is an addendum for the case of $s=-n\in\mathbb{Z}_{<0}$ giving the explicit form of the C\'{e}saro-remainder term under the action of $P^{n+1}$ in such cases - a formula and derivation we find combinatorially interesting in its own right.

In section 2.4 we then trivially extend our results to the setting of remainder summation relative to an arbitrary point $z_{0}\in\mathbb{C}$, rather than just the case of $z_{0}=0$ considered up to this point. This completes Part 1.

\section{C\'{e}saro-adapted scales and two theorems}

The two results from [1] which we will prove in this section are theorem 4 and theorem 3 from that paper. Writing $X=k+\alpha$ in the usual way from [1] (where our remainder sums for analysing the Riemann zeta function were all taken relative to the point $z_{0}=0$), these were as follows (relabelling them as theorems 1 and 2 respectively here):\\
\\
\textbf{Theorem 1:} \textit{For} $\delta\in\mathbb{C}$, $Re(\delta)\geq{0}$ \textit{and} $r\in\mathbb{Z_{\geq\textrm{0}}}$ \textit{we have that}:\\
\begin{equation}
k^{\delta}\alpha^{r}\,\overset{C}{\rightarrow}\begin{cases}
0 & \quad \textrm{if} \; \delta\notin\mathbb{Z_{\geq\textrm{0}}}\\
\frac{(-1)^{n}}{n+r+1} & \quad \textrm{if} \; \delta=n\in\mathbb{Z_{\geq\textrm{0}}}\end{cases}\label{CoreCesaroAsymptotic_k}\end{equation}\\ 
\\ \textit{Equivalently}\\
\begin{equation}X^{\delta}\alpha^{r}\,\overset{C}{\rightarrow}\begin{cases}
0 & \quad \textrm{if} \; \delta\neq{0}\\
\frac{1}{r+1} & \quad \textrm{if} \; \delta=0\end{cases}\label{CoreCesaroAsymptotic_X}\end{equation}\\
\\
\textbf{Theorem 2:} \textit{For any} $Re(s)\leq{1}, s\neq{1}$, \textit{we have locally uniformly the fundamental strong C\'{e}saro asymptotic relationship that}
\begin{equation}
s_{\zeta,s}(k+\alpha)\,=\,\sum_{j=1}^{k}j^{-s}\,\overset{C}{\simeq}\,\frac{X^{1-s}}{1-s}\,+\,\zeta(s)\,+\,o(1)\label{keystrongCesaroAsymptotic}\end{equation} \textit{via the pure power} $P^{r}$ \textit{where} $r=\lfloor Re(-s)\rfloor +1$.

\subsection{C\'{e}saro-adapted scales and proof of Theorem 1}

In [1] the initial result of theorem 1 (equation \ref{CoreCesaroAsymptotic_k}) was first derived for $k^{\delta}\alpha^{r}$ for $\delta=0,1,2,3,\ldots$ by using a process of iterative inversion at top-order\footnote{a la the theory of inversion of pseudo-differential operators}. It turns out, however, that it is easier to consider the case of $X^{\delta}\alpha^{r}$ first and then derive the result for $k^{\delta}\alpha^{r}$ from it.\footnote{this is an insight of Donald Cartwright} Deriving the result for $X^{\delta}\alpha^{r}$ is most cleanly done using the following concept:\\
\\
\textbf{Definition [C\'{e}saro-adapted scale]:} \textit{A C\'{e}saro-adapted scale is a 1-parameter family of functions} $\{q_{n}(X)\}_{n=1}^{\infty}$ \textit{such that}
\begin{equation}
\frac{d}{dX}q_{n+1}(X) = c_{n}\cdot q_{n}(X) \quad \textrm{for some set of non-zero constants} \; c_{n}, 
\label{Cesaro_adapted_scale_def_1}\end{equation}
\textit{and}
\begin{equation}
P[q_{n}](X) \rightarrow 0 \quad \textrm{classically as} \quad X\rightarrow\infty \quad \textrm{for all n}. 
\label{Cesaro_adapted_scale_def_2}\end{equation}
\\
To prove theorem 1, recall that when we took the p-sum function for $\zeta(0)$ in paper [1] we ended up re-writing it in the form
\begin{equation}
s(X)=s(k+\alpha)=\sum_{j=1}^{k}j^{0}=k=X+\zeta(0)-(\alpha-\frac{1}{2}) \quad .
\end{equation}
We thus take $q_{1}(X):=\alpha-\frac{1}{2}$ and build up a C\'{e}saro-adapted scale from it. Since functions defined in terms of $\alpha$ are automatically periodic with period $1$, the second requirement in the above definition that $P[q_{n}](X) \rightarrow 0$ classically as $X\rightarrow\infty$ means that we must have $\int_{0}^{1}q_{n}(\alpha)\,d\alpha\,=\,0$. Together with the fact that the first requirement in equation \ref{Cesaro_adapted_scale_def_1} means that each $q_{n+1}[X]$ is an anti-derivative of $q_{n}[X]$, this will uniquely define for us the constant of  integration at each stage in going from $q_{n}$ to $q_{n+1}$. We are thus able to calculate $q_{2}$ from $q_{1}$; then $q_{3}$ from $q_{2}$; and so on. We obtain
\begin{eqnarray*}
q_{1}(X) & = & \alpha-\frac{1}{2}\\
q_{2}(X) & = & \frac{1}{2}\alpha^{2}-\frac{1}{2}\alpha+\frac{1}{12}\\
q_{3}(X) & = & \frac{1}{3!}\alpha^{3}-\frac{1}{4}\alpha^{2}+\frac{1}{12}\alpha\\
q_{4}(X) & = & \frac{1}{4!}\alpha^{4}-\frac{1}{12}\alpha^{3}+\frac{1}{24}\alpha^{2}-\frac{1}{720}\\
 & \vdots &
\end{eqnarray*}
By construction these functions $\{q_{n}(X)\}_{n=1}^{\infty}$ form a C\'{e}saro-adapted scale with $c_{n}=1$ for all $n$ in equation \ref{Cesaro_adapted_scale_def_1}. Each $q_{n}(X)$ is a period-1 function $q_{n}(X)=q_{n}(k+\alpha)=q_{n}(\alpha)$ which is a polynomial of degree $n$ in $\alpha$, with leading term $\frac{1}{n!}\alpha^{n}$, and which is continuous for $n\geq{2}$, differentiable for $n\geq{3}$, twice-differentiable for $n\geq{4}$ and so on.\\
\\
\textbf{Proof of Theorem 1:} Using the $\{q_{n}(X)\}_{n=1}^{\infty}$, the proof that $X^{\delta}\alpha^{r}\rightarrow 0$ for $\delta\in\mathbb{C\setminus}\{0\}$ is then elementary. Taking $Re(\delta)\geq{0}$ without loss of generality (since trivially $X^{\delta}\alpha^{r} \rightarrow 0$ if $Re(\delta)<0$) we have that, for $\delta \neq{0}$,
\begin{eqnarray*}
P\left[\tilde{X}^{\delta}q_{r}(\tilde{X})\right](X) & = & \frac{1}{X}\int_{0}^{X}x^{\delta}q_{r}(x)\textrm{d}x \,=\, \frac{1}{X}\left\{\begin{array}{cc}\left[x^{\delta}q_{r+1}(x)\right]_{0}^{X} - \\
 \\
 \delta\int_{0}^{X}x^{\delta-1}q_{r+1}(x)\textrm{d}x\end{array}\right\}\\
 \\
 & = & X^{\delta-1}q_{r+1}(X)-\delta P\left[\tilde{X}^{\delta-1}q_{r+1}(\tilde{X})\right](X)\\
 \\
 & = & \left\{\begin{array}{cc}X^{\delta-1}q_{r+1}(X)-\delta X^{\delta-2}q_{r+2}(X) \\
 \\
 + \delta(\delta-1) P\left[\tilde{X}^{\delta-2}q_{r+2}(\tilde{X})\right](X)\end{array}\right\}\\
 \\
 & = & \left\{\begin{array}{cc}X^{\delta-1}q_{r+1}(X)-\delta X^{\delta-2}q_{r+2}(X) \\
 \\
 + \delta(\delta-1)X^{\delta-3}q_{r+3}(X)-\ldots \end{array}\right\}
\end{eqnarray*}
Since the highest power of $X$ in all the terms on the RHS is only $\delta-1$, it follows that by further repeated application of $P$ we can reduce $X^{\delta}q_{r}(X)$ to a sum of terms in which the power of $X$ has either negative or, at worst, zero real part. Specifically, if $\delta\notin\mathbb{Z}_{\geq{0}}$ then by taking $P^{\lfloor \delta \rfloor+1}\left[\tilde{X}^{\delta}q_{r}(\tilde{X})\right](X)$ we end up with an expression which is classically convergent to zero because it is a linear combination of terms either of the form $\frac{q_{r+j}(\alpha)}{X^{\rho}}$ with $Re(\rho)>0$ (which are clearly classically convergent to $0$); or of the images of such expressions under powers of $P$ (which remain classically convergent to $0$ since $P$ is regular). And if $\delta=n\in\mathbb{Z}_{\geq{0}}$ then on noting that $\{q_{j}(X)\}_{j=1}^{\infty}$ is a C\'{e}saro-adapted scale so that $P\left[q_{j}(\tilde{X})\right](X) \rightarrow 0$ classically for all $j$ as $X\rightarrow\infty$, it follows that $P^{n}\left[\tilde{X}^{n}q_{r}(\tilde{X})\right](X)=C \cdot q_{r+n}(X)+o(1)$ and so again $P^{n+1}\left[\tilde{X}^{n}q_{r}(\tilde{X})\right](X)$ is classically convergent to zero.

It follows that for arbitrary $Re(\delta)\geq{0}$ ($\delta\neq{0}$) and arbitrary $r\in\mathbb{Z}_{\geq{1}}$ we have
\begin{equation}
P^{\lfloor Re(\delta) \rfloor+1}\left[\tilde{X}^{\delta}q_{r}(\tilde{X})\right](X)=o(1) \quad \textrm{classically.}
\label{Cesaro_asympt_X_delta_q_r}\end{equation}
And when $r=0$, since $X^{\delta}$ is itself an eigenfunction of $P$ with eigenvalue $\frac{1}{\delta +1}$, so $X^{\delta}\overset{C}{\rightarrow}0$ under the regular polynomial $q_{\delta}(P)=\left(\frac{\delta+1}{\delta}\right)\left(P-\frac{1}{\delta+1}\right)$.

Result \ref{CoreCesaroAsymptotic_X} in theorem 1 therefore follows immediately on noting that clearly any power $\alpha^{r}$ can be expressed as a linear combination of the functions $\{q_{n}(\alpha)\}_{n=1}^{\infty}$ together with the constant function $1$. Here we are of course also noting that when $\delta=0$, it is trivial that $P[\alpha^{r}](X)\rightarrow \frac{1}{r+1}$ classically as $X\rightarrow\infty$.

But then, using result \ref{CoreCesaroAsymptotic_X}, it is straightforward to derive its discrete analogue given in equation \ref{CoreCesaroAsymptotic_k}. When $\delta=0$ this is simply the same formula as just noted, namely that $\alpha^{r}\overset{C}{\sim}\frac{1}{r+1}$. When $\delta\notin\mathbb{Z}_{\geq{1}}$, then writing $k^{\delta}\alpha^{r}$ as $(X-\alpha)^{\delta}\alpha^{r}=X^{\delta}\alpha^{r}-\binom{\delta}{1}X^{\delta-1}\alpha^{r+1}+\binom{\delta}{2}X^{\delta-2}\alpha^{r+2}-\ldots$ gives $k^{\delta}\alpha^{r}$ as a sum of a finite collection of classically divergent terms plus a classically $o(1)$ residual; and invocation of result \ref{CoreCesaroAsymptotic_X} for each of these divergent terms $X^{\delta}\alpha^{r}$, $X^{\delta-1}\alpha^{r+1}, \ldots$ in turn immediately yields that $k^{\delta}\alpha^{r}\overset{C}{\rightarrow}0$ as claimed. As for the case of $\delta=n\in\mathbb{Z}_{\geq{0}}$, this finally follows in the same way on noting that $k^{n}=(X-\alpha)^{n}$ is now a finite expression in which the only term with no residual power of $X$ is $(-1)^{n}\alpha^{n}$ and then invoking result \ref{CoreCesaroAsymptotic_X} to conclude that
\begin{equation*}
\underset{X\rightarrow\infty}{Clim}\,k^{\delta}\alpha^{r} = \underset{X\rightarrow\infty}{Clim}\,(-1)^{n}\alpha^{n+r} = \frac{(-1)^{n}}{n+r+1} \quad .
\end{equation*}
This completes the proof of theorem 1.\\
\\
\textbf{Comments:} \textbf{(i)} Note that, in the course of the above derivation, we have in fact proved that the C\'{e}saro asymptotic result in equation \ref{CoreCesaroAsymptotic_X} arises via the regular polynomial $q_{\delta}(P):=\left(\frac{\delta+1}{\delta}\right)\left(P-\frac{1}{\delta+1}\right)\cdot P^{\lfloor Re(\delta) \rfloor + 1}$. This is because, while equation \ref{Cesaro_asympt_X_delta_q_r} is a strong C\'{e}saro  relationship, re-expressing $\alpha^{r}$ in terms of the basis consisting of the constant function $1$ and the $\{q_{n}(\alpha)\}_{n=1}^{\infty}$ leads to a stand-alone term $X^{\delta}$ whose C\'{e}saro  annihilation then requires the factor $\left(\frac{\delta+1}{\delta}\right)\left(P-\frac{1}{\delta+1}\right)$. As for equation \ref{CoreCesaroAsymptotic_k}, since its proof involved consideration of a finite set of expressions $X^{\delta-j}\alpha^{r+j}$, $j=0,1, \ldots ,\lfloor \delta \rfloor$, each of which may lead to a stand-alone term $X^{\delta}$, $X^{\delta-1}$, $X^{\delta-2}$ etc, so the regular polynomial giving this C\'{e}saro asymptotic result involves additional factors of the form $\left(\frac{\delta-j+1}{\delta-j}\right)\left(P-\frac{1}{\delta-j+1}\right)$.\\
\\
\textbf{(ii)} The key in the above proof was first to focus on the continuous case of $X^{\delta}\alpha^{r}$ before tackling the discrete case of $k^{\delta}\alpha^{r}$; and secondly to change from considering pure powers $1,\alpha,\alpha^{2},\alpha^{3},\ldots$ in such expressions to considering an alternative basis of polynomials in $\alpha$, namely the C\'{e}saro-adapted set of functions $1,q_{1}(\alpha),q_{2}(\alpha),q_{3}(\alpha),\ldots$. Before turning to the proof of theorem 2, we take a moment to consider these functions in greater depth.

\subsection{The functions $q_{n}(X)$, $\tilde{q}_{n}(X)$, $\overset{\lor}{q}_{n}(X)$; Bernoulli polynomials; Formal symbols and formal function elements}

\textbf{The functions $q_{n}(X)$, $\tilde{q}_{n}(X)$, $\overset{\lor}{q}_{n}(X)$:} The set of functions $\{q_{n}(X)\}_{n=1}^{\infty}$ can be normalised in various alternate ways while remaining a C\'{e}saro-adapted set, since doing so just affects the values of the constants, $c_{n}$, appearing in equation \ref{Cesaro_adapted_scale_def_1}. In particular it will be convenient to consider the set $\{\tilde{q}_{n}(X)\}_{n=1}^{\infty}$ given by
\begin{equation}
\tilde{q}_{n}(X):=(n-1)!\,q_{n}(X)=(n-1)!\,q_{n}(\alpha)
\label{q_tilde_def}\end{equation}
so that $\tilde{q}_{1}(\alpha)=\alpha-\frac{1}{2}$, $\tilde{q}_{2}(\alpha)=\frac{1}{2}\alpha^{2}-\frac{1}{2}\alpha+\frac{1}{12}$, $\tilde{q}_{3}(\alpha)=\frac{1}{3}\alpha^{3}-\frac{1}{2}\alpha^{2}+\frac{1}{6}\alpha$, $\tilde{q}_{4}(\alpha)=\frac{1}{4}\alpha^{4}-\frac{1}{2}\alpha^{3}+\frac{1}{4}\alpha^{2}-\frac{1}{120}$ and so on. The polynomials $\tilde{q}_{n}(\alpha)$ always have leading terms $\frac{1}{n}\alpha^{n}-\frac{1}{2}\alpha^{n-1}$ and constant term $(-1)^{n-1}\zeta(1-n)$.

It will also be of interest to define a "shifted" version of these polynomials, namely the set $\{\overset{\lor}{q}_{n}(X)\}_{n=0}^{\infty}$ given by
\begin{equation}
\overset{\lor}{q}_{n}(X):=\tilde{q}_{n+1}(X) \, .
\label{q_cech_def}\end{equation}
The $\{\overset{\lor}{q}_{n}\}$ seem a curious set to introduce since they are merely the same as the set $\{\tilde{q}_{n}\}$ but their polynomial degrees are now dislocated so that $\overset{\lor}{q}_{n}(\alpha)$ no longer has degree $n$, but rather $n+1$. We will shortly see, however, that there is an alternative perspective under which the $\{\overset{\lor}{q}_{n}\}$ become the most natural formulation of this set of functions.\\
\\
\textbf{Bernoulli numbers and polynomials:} Recall also the definition of the well-known Bernoulli numbers and Bernoulli polynomials. The Bernoulli numbers are given recursively by
\begin{equation}
B_{0}=1\quad\textrm{and}\quad B_{n}=\sum_{j=0}^{n}\dbinom{n}{j}\, B_{j}\quad\forall n\geq2
\label{eq:BernoulliNums1}\end{equation}
which implies $B_{1}=\frac{-1}{2}$, $B_{2}=\frac{1}{6}$, $B_{4}=\frac{-1}{30}$, $B_{6}=\frac{1}{42}$ and so on. For $n\geq1$ we have $B_{2n+1}=0$ and the even-index Bernoulli numbers are related to the values of $\zeta$ at negative odd integers by $B_{2n}=-2n \cdot \zeta(-2n+1)$. A generating function for the $B_{n}$ is given by
\begin{equation}
\frac{t}{\textrm{e}^{t}-1}=\sum_{n=0}^{\infty}\frac{B_{n}}{n!}\, t^{n} \; .\label{eq:BernoulliNumsGenFn}\end{equation}
Bernoulli polynomials are then defined by
\begin{equation}
B_{n}(x):=\sum_{j=0}^{n}\dbinom{n}{j}\, B_{j}x^{n-j}=\sum_{j=0}^{n}\dbinom{n}{j}\, B_{n-j}x^{j}
\label{eq:BernoulliPolys1}\end{equation}
and are monic polynomials satisfying $B_{n}(0)=B_{n}=B_{n}(1)$. Since $B_{n}(0)=B_{n}(1)$ we alternatively let $\tilde{B}_{n}(x)$ be the periodic, period-1 extensions of the Bernoulli polynomials from $[0,1]$ to all of $\mathbb{R}$. For $n\geq3$ these are differentiable and satisfy 
$\frac{\textrm{d}}{\textrm{d}x}\tilde{B}_{n}(x)=n\tilde{B}_{n-1}(x)$
and this relationship continues to hold for $n=2$ except at the integer points (of measure $0$) where $\tilde{B}_{1}(x)$ is defined to have value $0$. The first few $B_{n}(x)$ and $\tilde{B}_{n}(x)$ are given by
\begin{eqnarray}
B_{1}(x)=x-\frac{1}{2} & , & \tilde{B}_{1}(x)=\alpha-\frac{1}{2}\nonumber \\
\nonumber \\B_{2}(x)=x^{2}-x+\frac{1}{6} & , & \tilde{B}_{2}(x)=\alpha^{2}-\alpha+\frac{1}{6}\nonumber \\
\nonumber \\B_{3}(x)=x^{3}-\frac{3}{2}x^{2}+\frac{1}{2}x & , & \tilde{B}_{3}(x)=\alpha^{3}-\frac{3}{2}\alpha^{2}+\frac{1}{2}\alpha
\label{eq:BernoulliPolys2}\end{eqnarray}
and so on where $\alpha=x-\lfloor x \rfloor = x-k$ is the usual fractional part of $x$. We denote by $b_{n}(k)$ the polynomials in $k\in\mathbb{Z_{>\textrm{0}}}$ given by the sums of $(n-1)^{st}$ powers:
\begin{equation}
b_{n}(k):=\sum_{j=1}^{k}\, j^{n-1}
\label{eq:SumPowersPolys}\end{equation}
which are closely related to the Bernoulli polynomials.

We see that our function sets $\{q_{n}\}$, $\{\tilde{q}_{n}\}$ and $\{\overset{\lor}{q}_{n}\}$ are all just alternative normalisations of the periodised Bernoulli polynomials; or alternatively that a new perspective regarding the periodised Bernoulli polynomials is that they are a natural, periodic C\'{e}saro-adapted set, and one which is particularly well-fitted to the C\'{e}saro analysis of the zeta function using theorems 1 and 2 (see [1]).\\
\\
\textbf{Formal symbols and formal function elements:} Connecting to the Bernoulli numbers and polynomials leads in one further interesting direction. If we introduce a formal symbol, $B$, satisfying
\begin{equation}
B^{n}:=B_{n}
\label{B_formal_symbol_def_1}\end{equation}
then equation \ref{eq:BernoulliNums1} becomes simply
\begin{equation}
B^{n}:=(1+B)^{n} \qquad \textrm{for} \quad n\in\mathbb{Z}_{\geq{2}} \quad .
\label{eq:BernoulliNums2}\end{equation}
Thus the defining relation for the Bernoulli numbers takes on a particularly compact form when expressed in terms of $B$. The symbol $B$ is of course purely formal - the Bernoulli numbers are certainly not powers of any actual complex number - but the use of $B$ simplifies their defining relation and so seems the natural way of defining them.

In the same spirit, we introduce two further formal symbols which are closely related. Let the symbol $\zeta$ satisfy
\begin{equation}
\zeta^{\rho}:=\zeta(\rho) \quad \textrm{for all} \quad \rho\in\mathbb{C\setminus}\{1\}
\label{zeta_formal_symbol_def_1}\end{equation}
and let $\tau$ be its formal inverse given by
\begin{equation}
\tau^{\rho}:=\zeta(-\rho) \quad \textrm{for all} \quad \rho\in\mathbb{C\setminus}\{-1\} \, .
\label{zeta_formal_symbol_def_1}\end{equation}
Using $\tau$ we can express the polynomials $\{\overset{\lor}{q}_{n}(\alpha)\}_{n=0}^{\infty}$ in a very simple form as long as we are careful in our interpretation of the binomial theorem. 

Consider the binomial expression $(\alpha-\tau)^{\rho}$ mixing $\alpha \in [0,1)$ and the formal symbol $\tau$. Since we always have $\alpha$ small we consider its expansion as a Taylor series in non-negative powers of $\alpha$:
\begin{equation}
(\alpha-\tau)^{\rho}\,=\,\sum_{j=0}^{\infty}\binom{\rho}{j}\alpha^{j}(-1)^{\rho-j}\tau^{\rho-j}\,=\,\sum_{j=0}^{\infty}e^{i\pi (\rho-j)}\binom{\rho}{j}\zeta(j-\rho)\alpha^{j} \, .
\label{alpha_minus_tau_to_rho_1}\end{equation}

Now it is well-known (See appendix 1) that for large $j$, $\binom{\rho}{j}$ has an asymptotic expansion in descending powers $j^{-\rho-1}$, $j^{-\rho-2}$, $j^{-\rho-3}$ etc:
\begin{equation*}
\binom{\rho}{j}\,=\,a_{1}(\rho)j^{-\rho-1}+a_{2}(\rho)j^{-\rho-2}+a_{3}(\rho)j^{-\rho-3}+ \ldots \, .
\end{equation*}
Since $\zeta(j-\rho)\rightarrow 1$ as $j\rightarrow\infty$ it follows that the RHS in equation \ref{alpha_minus_tau_to_rho_1} will give a well-defined, classically convergent Taylor series with radius of convergence $1$ for $Re(\rho)>0$. 

And it is equally clear that, if we use generalised C\'{e}saro convergence via regular polynomials $q_{\rho}(P)$ in the same manner as we used it throughout the first two papers in this introductory series ([1] and [2]), we can then analytically continue to $Re(\rho)\leq{0}$ one strip at a time to give a well-defined meaning to $(\alpha-\tau)^{\rho}$ at least for all $\rho\in\mathbb{C\setminus}\mathbb{Z}_{\leq{0}}$. For example, to extend to the strip $-1<Re(\rho)\leq{0}$ we would use $q_{\rho}(P)=\left(\frac{\rho}{\rho+1}\right)\left(P+\frac{1}{\rho}\right)$; and then for the strip $-2<Re(\rho)\leq{-1}$ we would use $q_{\rho}(P)=\left(\frac{\rho}{\rho+2}\right)\left(P+\frac{1}{\rho}\right)\left(P+\frac{1}{\rho+1}\right)$; and so on.

We will return to consider $(\alpha-\tau)^{\rho}$ as a general function of the complex variable $\rho$ in section 4, but for now let us focus just on the cases $\rho=n\in\mathbb{Z}_{\geq{0}}$.

When $\rho=0$, equation \ref{alpha_minus_tau_to_rho_1} gives us formally
\begin{equation*}
(\alpha-\tau)^{0}\,=\,\binom{0}{0}\tau^{0}\alpha^{0}-\binom{0}{1}\tau^{-1}\alpha^{1}+\binom{0}{2}\tau^{-2}\alpha^{2}- \ldots
\end{equation*}
and since for $j\geq{2}$ we have that $\binom{0}{j}=0$ and $\tau^{-j}=\zeta(j)$ is finite, this reduces to just the first two terms
\begin{equation*}
(\alpha-\tau)^{0}\,=\,\binom{0}{0}\tau^{0}\alpha^{0}-\binom{0}{1}\tau^{-1}\alpha^{1} \, .
\end{equation*}
The first of these is a constant term $\tau^{0}=\zeta(0)=-\frac{1}{2}$ but for the second we need to interpret the coefficient of $\alpha$ as $-\lim_{\epsilon\rightarrow0}\binom{0}{1-\epsilon}\tau^{-(1-\epsilon)}$ in order to make sense of the formal product of the zero from $\binom{0}{1}$ and the pole from $\tau^{-1}=\zeta(1)$. We get $-\lim_{\epsilon\rightarrow0}\frac{0!}{(1-\epsilon)!(-1+\epsilon)!} \zeta(1-\epsilon)=-\lim_{\epsilon\rightarrow0} \epsilon\cdot\left(\frac{-1}{\epsilon}\right)=1$ and so the coefficient of $\alpha$ is $1$. We thus get
\begin{equation*}
(\alpha-\tau)^{0} \,=\, \alpha - \frac{1}{2} \,=\, \overset{\lor}{q}_{0}(\alpha) \quad .
\end{equation*}
In the same fashion, using the same limit-based interpretation to make sense of the coefficient of the leading term, we get
\begin{eqnarray*}
(\alpha-\tau)^{1} & = & -\binom{1}{0}\tau^{1}\alpha^{0}+\binom{1}{1}\tau^{0}\alpha^{1}-\binom{1}{2}\tau^{-1}\alpha^{2}\\
 & = & \frac{1}{2}\alpha^{2}-\frac{1}{2}\alpha+\frac{1}{12} \,=\, \overset{\lor}{q}_{1}(\alpha)
\end{eqnarray*}
and likewise
\begin{eqnarray*}
(\alpha-\tau)^{2} & = & \frac{1}{3}\alpha^{3}-\frac{1}{2}\alpha^{2}+\frac{1}{6}\alpha \,=\, \overset{\lor}{q}_{2}(\alpha) \, ,\\
(\alpha-\tau)^{3} & = & \frac{1}{4}\alpha^{4}-\frac{1}{2}\alpha^{3}+\frac{1}{4}\alpha^{2}-\frac{1}{120} \,=\, \overset{\lor}{q}_{3}(\alpha)
\end{eqnarray*}
and in general, for all $n\in\mathbb{Z}_{\geq{0}}$,
\begin{equation}
(\alpha-\tau)^{n}\,=\,\overset{\lor}{q}_{n}(\alpha)\,=\,\tilde{q}_{n+1}(\alpha) \quad .
\label{alpha_minus_tau_to_n_1}\end{equation}

This equation provides the sense, noted previously, in which the apparently oddly-defined set $\{\overset{\lor}{q}_{n}\}$ are in fact the most naturally-defined set of functions among our various collections of differently-normalised C\'{e}saro-adapted sets. After introducing the formal symbol $\tau$ they have the simplest and cleanest form, namely $\overset{\lor}{q}_{n}(\alpha)=(\alpha-\tau)^{n}$.

We shall return to undertake further investigations regarding formal symbols like $\tau$, and how they seem to provide the "right" framework for interpreting and expressing many of our results, in section 4.

For now, however, we wish to return our focus towards proving theorem 2, and to that end we introduce one further formal extension:\\
\\
\textbf{Definition 1 [Formal function elements]:} \textit{Just as} $B$, $\zeta$ \textit{and} $\tau$ \textit{are formal} \textrm{symbols} \textit{which resolve into complex numbers when evaluated after powering, we also define formal} \textrm{function elements} $q$, $\tilde{q}$, $\overset{\lor}{q}$ \textit{and} $\tilde{B}$ \textit{which resolve into functions when evaluated after powering, i.e.}
\begin{equation}
\tilde{q}^{n}\,:=\,\tilde{q}_{n}
\label{Function_elements_def_1}\end{equation}
\textit{and similarly for} $q$, $\overset{\lor}{q}$ \textit{and} $\tilde{B}$.\\
\\
\textbf{Comment:} Given that we have seen that $\overset{\lor}{q}_{n}(\alpha)=(\alpha-\tau)^{n}=\tilde{q}_{n+1}(\alpha)$ and that $(\alpha-\tau)^{\rho}$ makes sense as a well-defined function of $\alpha$ for $\rho\in\mathbb{C\setminus}\mathbb{Z}_{<0}$, we in fact generalise our definition of $\tilde{q}$ and $\overset{\lor}{q}$ to have
\begin{equation}
\overset{\lor}{q}_{\rho}(\alpha):=(\alpha-\tau)^{\rho}=\tilde{q}_{\rho+1}(\alpha)
\label{Function_elements_def_1_rho}\end{equation}
and then extend the function element definition in equation \ref{Function_elements_def_1} to have
\begin{equation}
\tilde{q}^{\rho}\,:=\,\tilde{q}_{\rho} \quad \textrm{for all $\rho\in\mathbb{C\setminus}\mathbb{Z}_{<0}$}
\label{Function_elements_def_2}\end{equation}
and similarly for $\overset{\lor}{q}$, $q$ and $\tilde{B}$.

What we need in order to prove theorem 2, however, is just the case of $n\in\mathbb{Z}_{\geq{0}}$ in equation \ref{Function_elements_def_1} and so we defer following these intriguing paths into the formal jungles until Part 2.

\subsection{Proof of Theorem 2}

The expression which we have used in [1] and [2] for the p-sum of the defining series for $\zeta(s)$ is just its associated Euler-McLaurin asymptotic series in k:
\begin{equation}
s_{\zeta,s}(k+\alpha)\,=\,\frac{k^{1-s}}{1-s}\,+\,\zeta(s)\,+\,\frac{1}{2}k^{-s}-\,\sum_{r=1}^{\infty}\frac{B_{2r}}{(2r)!}s(s+1)\dots(s+2r-2)k^{-s-2r+1}
\label{psumzeta_1}\end{equation}
The key to proving theorem 2 is to recast this into an an exact expression in terms of $X$ and our C\'{e}saro-adapted set $\{\overset{\lor}{q}_{n}\}$. Writing $X=k+\alpha$ as always, and replacing $k$ everywhere by $X-\alpha$ we obtain
\begin{eqnarray*}
s_{\zeta,s}(k+\alpha) & = & \zeta(s)\,+\,\frac{1}{1-s} \left\{\begin{array}{cc}X^{1-s}-(1-s)X^{-s}\alpha+\frac{(1-s)(-s)}{2!}X^{-s-1}\alpha^{2}\\
 \\
 -\frac{(1-s)(-s)(-s-1)}{3!}X^{-s-2}\alpha^{3}+\ldots \end{array}\right\} \\
 \\
 &   & +\frac{1}{2}\left\{X^{-s}-(-s)X^{-s-1}\alpha+\frac{(-s)(-s-1)}{2!}X^{-s-2}\alpha^{2}-\ldots \right\} \\
 \\
 &   & -\frac{s}{12}\left\{X^{-s-1}-(-s-1)X^{-s-2}\alpha+\ldots \right\} \\
 \\
 &   & +\frac{s(s+1)(s+2)}{720}\left\{X^{-s-3}-\ldots \right\} - \ldots \\
 \\
 & = & \frac{X^{1-s}}{1-s} + \zeta(s) - X^{-s}\left(\alpha-\frac{1}{2}\right) - s X^{-s-1}\left(\frac{1}{2}\alpha^{2}-\frac{1}{2}\alpha+\frac{1}{12}\right)\\
 &   & - \frac{s(s+1)}{2!} X^{-s-2}\left(\frac{1}{3}\alpha^{3}-\frac{1}{2}\alpha^{2}+\frac{1}{6}\alpha\right) + \ldots
 \end{eqnarray*}
\begin{equation}
= \frac{X^{1-s}}{1-s} + \zeta(s) - \left\{X^{-s}\overset{\lor}{q}_{0}(\alpha) - \binom{-s}{1} X^{-s-1}\overset{\lor}{q}_{1}(\alpha) + \binom{-s}{2} X^{-s-2}\overset{\lor}{q}_{2}(\alpha) - \ldots \right\}
\label{psumzeta_X_q_cech_exact_1}\end{equation}
which can be very compactly expressed in terms of the function element $\overset{\lor}{q}$ as
\begin{equation}
s_{\zeta,s}(k+\alpha) = \frac{X^{1-s}}{1-s} + \zeta(s) - \left(X-\overset{\lor}{q}\right)^{-s} \, .
\label{psumzeta_X_q_cech_exact_2}\end{equation}
Now in the course of proving theorem 1 (the argument is unchanged by the use of the C\'{e}saro-adapted set $\{\overset{\lor}{q}_{n}\}$ rather than the set $\{q_{n}\}$) we saw that for any $\delta\in\mathbb{C\setminus}\{0\}$ and any $r\in\mathbb{Z}_{\geq{0}}$ we have
\begin{equation}
P^{\lfloor Re(\delta) \rfloor + 1}\left[\tilde{X}^{\delta}\overset{\lor}{q}_{r}(\tilde{\alpha})\right] \rightarrow 0 \quad \textrm{classically as} \quad X\rightarrow\infty \quad .
\label{X_delta_q_cech_r_lemma}\end{equation}
It follows that we have the strong C\'{e}saro asymptotic relationship
\begin{equation*}
s_{\zeta,s}(k+\alpha) = \sum_{j=1}^{k}j^{-s} \, \overset{C}{\simeq} \, \frac{X^{1-s}}{1-s} + \zeta(s) + o(1)
\end{equation*}
via the pure power $P^{\lfloor Re(-s) \rfloor + 1}$. This completes the proof of theorem 2. \\
\\
\textbf{Comments: (i)} This proof also confirms the claim from [1] that the regular polynomial $q(P;s)$ which provides the C\'{e}saro convergence of the p-sum function $s_{\zeta,s}(X) = \sum_{j=1}^{k}j^{-s}$ to $\zeta(s)$ is given by $q(P;s)=\left(\frac{2-s}{1-s}\right)\left(P-\frac{1}{2-s}\right)P^{\lfloor Re(-s) \rfloor + 1}$ for all $s\in\mathbb{C\setminus}\{1\}$.\\
\\
\textbf{(ii)} When $Re(s)\leq{1}$ but $s\notin \mathbb{Z}_{\leq{1}}$ then the expression $\left(X-\overset{\lor}{q}\right)^{-s}=X^{-s}\overset{\lor}{q}_{0}(\alpha)+sX^{-s-1}\overset{\lor}{q}_{1}(\alpha)+\ldots \,$ has infinitely many terms, but only finitely many where the power of $X$ has positive real part, after which all terms involve a power of $X$ with negative real part and are therefore already classically convergent to $0$.

When $s=-n\in\mathbb{Z}_{\leq{0}}$, however, the expression $\left(X-\overset{\lor}{q}\right)^{n}=X^{n}\overset{\lor}{q}_{0}(\alpha)-nX^{n-1}\overset{\lor}{q}_{1}(\alpha)+\ldots+(-1)^{n}X^{0}\overset{\lor}{q}_{n}(\alpha) \,$ has only a finite number of terms, and it is interesting to see what $P^{n+1}\left[\left(\tilde{X}-\overset{\lor}{q}\right)^{n}\right](X)$ produces in these circumstances as an exact expression, beyond just its asymptotic convergence to $0$ as $X\rightarrow\infty$. 

We omit details of the derivations, but the final formulae obtained are again compact and clean. We have, first, that
\begin{eqnarray}
P^{n}\left[\tilde{X}^{n}\overset{\lor}{q}_{r}(\tilde{\alpha})\right] & = & \frac{(-1)^{n}n!}{(r+1)\cdots(r+n)}(P-1)(P-\frac{1}{2})\cdots(P-\frac{1}{n})\left[\overset{\lor}{q}_{r+n}(\tilde{\alpha})\right] \nonumber\\
 & = & \frac{(-1)^{n}}{\binom{r+n}{n}}(P-1)(P-\frac{1}{2})\cdots(P-\frac{1}{n})\left[\overset{\lor}{q}_{r+n}(\tilde{\alpha})\right] \quad .
\label{X_n_q_cech_r_calc_1}\end{eqnarray}
From this it follows that
\begin{equation}
P^{n}\left[\left(\tilde{X}-\overset{\lor}{q}\right)^{n}\right]=(-1)^{n}(n+1)(P-\frac{1}{2})(P-\frac{1}{3})\cdots(P-\frac{1}{n+1})\left[\overset{\lor}{q}_{n}(\tilde{\alpha})\right]
\label{X_minus_q_cech_to_n_calc_1}\end{equation}
on noting the purely combinatorial identity that
\begin{equation}
\sum_{j=0}^{n}(P-1)(P-\frac{1}{2})\cdots(P-\frac{1}{n-j})\cdot P^{j} = (n+1)(P-\frac{1}{2})(P-\frac{1}{3})\cdots(P-\frac{1}{n+1}) \, .
\label{X_minus_q_cech_to_n_calc_2}\end{equation}
Hence at once from equation \ref{X_minus_q_cech_to_n_calc_1} we have the exact expression
\begin{equation}
P^{n+1}\left[\left(\tilde{X}-\overset{\lor}{q}\right)^{n}\right]=(-1)^{n}(n+1)(P-\frac{1}{2})(P-\frac{1}{3})\cdots(P-\frac{1}{n+1})\cdot P\left[\overset{\lor}{q}_{n+1}(\tilde{\alpha})\right]
\label{X_minus_q_cech_to_n_calc_3}\end{equation}
and the classical convergence to $0$ of $\left(X-\overset{\lor}{q}\right)^{n}$ under $P^{n+1}$ is immediate from this since $P\left[\overset{\lor}{q}_{n+1}(\tilde{\alpha})\right] \rightarrow 0$ as $X\rightarrow\infty$ and $(n+1)(P-\frac{1}{2})(P-\frac{1}{3})\cdots(P-\frac{1}{n+1})$ is a regular operator.

We find the formulae \ref{X_n_q_cech_r_calc_1} - \ref{X_minus_q_cech_to_n_calc_2} intriguing and it would be interesting to explore how analogous formulae might be developed for the case of $s\notin\mathbb{Z}_{<0}$.\\
\\
\textbf{(iii)} The exact expression given in equation \ref{psumzeta_X_q_cech_exact_2} is not just central to proving theorem 2 and hence to our analysis of $\zeta$ and $\zeta_{H}$ in the first two of these introductory papers. It will also be the key building block in developing the concept of C\'{e}saro arrays, which will be the focus of our next set of three papers.\\
\\
\subsection{Extending theorems 1 and 2 to the remainder C\'{e}saro setting}

Theorems 1 and 2 both relate to p-sum functions defined on the positive real axis, so that $X=k+\alpha$ lies on $\mathbb{R}_{\geq{0}}$. The more general setting for C\'{e}saro analysis, however, consists of \textit{remainder} summation relative to an arbitrary point $z_{0}\in\mathbb{C}$, with continuous summation variable $z=z_{0}+k+\alpha$ lying on a horizontal ray in the complex plane starting at $z_{0}$ and going to $\infty$. Viewed this way, these theorems both pertain just to the case of $z_{0}=0$ and it is worth considering how they both generalize to arbitrary $z_{0}$.

In fact, such a generalization is relatively straightforward and the analysis involved closely follows that undertaken in sections 2.1-2.3. If we fix $z_{0}$ and let $z=z_{0}+k+\alpha$ then, for theorem 1, the working to derive equation \ref{Cesaro_asympt_X_delta_q_r} is easily adapted and we find that for arbitrary $Re(\delta)\geq{0}$ and $r\in\mathbb{Z}_{\geq{1}}$ we still have
\begin{equation}
P^{\lfloor Re(\delta) \rfloor+1}\left[\tilde{z}^{\delta}q_{r}(\tilde{\alpha})\right](z)\rightarrow 0 \quad \textrm{classically as} \; z\rightarrow\infty \, .
\label{Cesaro_asympt_X_delta_q_r_2}\end{equation}
It follows as before that, for any $\delta\in\mathbb{C\setminus}\{0\}$ we have $z^{\delta}\alpha^{r}\overset{C}{\rightarrow}0$; while, for $\delta=0$, $\alpha^{r}\overset{C}{\rightarrow}\frac{1}{r+1}$ in the usual way. Thus the analogue of equation \ref{CoreCesaroAsymptotic_X} in theorem 1 holds in identical fashion in this more general context.

As for the analogue of equation \ref{CoreCesaroAsymptotic_k} for $k^{\delta}\alpha^{r}$ in this context, we now write $k$ as $z-(z_{0}+\alpha)$ so that we have
\begin{equation}
k^{\delta}\alpha^{r}=z^{\delta}\alpha^{r}-\binom{\delta}{1}z^{\delta-1}(z_{0}+\alpha)\alpha^{r}+\binom{\delta}{2}z^{\delta-2}(z_{0}+\alpha)^{2}\alpha^{r}-\ldots
\label{k_delta_alpha_r_remainder_expansion}\end{equation}
For $\delta\notin\mathbb{Z}_{\geq{0}}$ it follows immediately from the result just derived for $z^{\delta}\alpha^{r}$ that we therefore still have $k^{\delta}\alpha^{r}\overset{C}{\rightarrow}0$. As for $\delta=n\in\mathbb{Z}_{\geq{0}}$, in this case equation \ref{k_delta_alpha_r_remainder_expansion} becomes a finite expansion and it then follows that the only term with no residual positive power of $z$ is the term $(-1)^{n}\binom{n}{n}(z_{0}+\alpha)^{n}\alpha^{r}$; and since we still have trivially that $P\left[\tilde{\alpha}^{m}\right](z)\rightarrow\frac{1}{m+1}$, it follows that $k^{n}\alpha^{r}\overset{C}{\rightarrow}(-1)^{n}\sum_{j=0}^{n}\binom{n}{j}z_{0}^{n-j}\cdot \frac{1}{r+j+1}$.

Overall, we have derived the following general analogue of theorem 1 for remainder C\'{e}saro summation relative to $z_{0}$:\\
\\
\textbf{Theorem 1a:} \textit{If we fix} $z_{0}\in\mathbb{C}$ \textit{and let} $z=z_{0}+k+\alpha$ \textit{then for} $\delta\in\mathbb{C}$, $Re(\delta)\geq{0}$ \textit{and} $r\in\mathbb{Z_{\geq\textrm{0}}}$ \textit{we have that}:\\
\begin{equation}z^{\delta}\alpha^{r}\,\overset{C}{\rightarrow}\begin{cases}
0 & \quad \textrm{if} \; \delta\neq{0}\\
\frac{1}{r+1} & \quad \textrm{if} \; \delta=0\end{cases}\label{CoreCesaroAsymptotic_X_remainder_setting}\end{equation}\\ 
\textit{Equivalently}\\
\begin{equation}
k^{\delta}\alpha^{r}\,\overset{C}{\rightarrow}\begin{cases}
0 & \quad \textrm{if} \; \delta\notin\mathbb{Z_{\geq\textrm{0}}}\\
(-1)^{n}\sum_{j=0}^{n}\binom{n}{j}z_{0}^{n-j}\cdot \frac{1}{r+j+1} & \quad \textrm{if} \; \delta=n\in\mathbb{Z_{\geq\textrm{0}}}\end{cases}\label{CoreCesaroAsymptotic_k_remainder_setting}\end{equation}\\
\\
As for theorem 2 concerning the p-sum of the Riemann zeta function, $\zeta(s)$, the analogue in this case clearly concerns the p-sum of the Hurewicz zeta function, $\zeta_{H}(z_{0};s)$, at an arbitrary fixed $z_{0}\in\mathbb{C}$. In this case the working leading to equations \ref{psumzeta_X_q_cech_exact_1} and \ref{psumzeta_X_q_cech_exact_2} applies in identical fashion to give us that
\begin{equation}
s_{\zeta_{H},s}(z)=\sum_{j=1}^{k}(z_{0}+j)^{-s}=\zeta_{H}(z_{0};s)+\frac{z^{1-s}}{1-s}+\left(z-\overset{\lor}{q}\right)^{-s} \; .
\label{psumzeta_z_q_cech_exact_1}\end{equation}
It follows from result \ref{Cesaro_asympt_X_delta_q_r_2} that an exact analogue of theorem 2 holds in this remainder C\'{e}saro context:\\
\\
\textbf{Theorem 2a:} \textit{For any given} $z_{0}\in\mathbb{C}$ \textit{and for any} $Re(s)\leq{1}, s\neq{1}$, \textit{we have locally uniformly the fundamental strong C\'{e}saro asymptotic relationship that}
\begin{equation}
s_{\zeta_{H},s}(z)\,=\,\sum_{j=1}^{k}(z_{0}+j)^{-s}\,\overset{C}{\simeq}\,\frac{z^{1-s}}{1-s}\,+\,\zeta_{H}(z_{0};s)\,+\,o(1)\label{keystrongCesaroAsymptotic}\end{equation} \textit{via the pure power} $P^{r}$ \textit{where} $r=\lfloor Re(-s)\rfloor +1$.\\
\\
In identical fashion to the logic applied to obtain $\zeta(s)$ as the C\'{e}saro limit of the p-sum $\sum_{j=1}^{k}j^{-s}$, it of course follows immediately from theorem 2a that the Hurewicz zeta function, $\zeta_{H}(z_{0};s)$, is given by the C\'{e}saro limit of the p-sum $\sum_{j=1}^{k}(z_{0}+j)^{-s}$ via the same regular polynomial, $q(P;s)=\left(\frac{2-s}{1-s}\right)\left(P-\frac{1}{2-s}\right)\cdot P^{\lfloor Re(-s) \rfloor + 1}$. Note that this polynomial is independent of $z_{0}$.

\section{Introduction to Part 2: The formal jungles}

The primary purpose of this paper was to introduce the key concept of a C\'{e}saro-adapted scale and to use it to prove theorems 1 and 2 concerning generalised C\'{e}saro analysis. In particular, the strong C\'{e}saro asymptotic result of theorem 2 will be the key building block in developing the notion of C\'{e}saro arrays in our next papers. This work is done and the reader is entitled to stop at this point and rest from his labours if he wishes.

At the same time, we have seen in section 2 how much the introduction of the ideas of \textit{formal symbols} and \textit{formal function elements} has simplified the expression of key formulae and arguments which are otherwise somewhat messy. In a sense the use of such concepts seems to be "the right way" to view these results; and we have caught glimpses of how useful they can be.

As such, we conclude this last of our introductory papers with a brief "aerial tour" of the formal jungles - partly in order to explore these concepts and to note further interesting results and simplifications; partly to extend our understanding of the family of functions $\overset{\lor}{q}_{\rho}(\alpha)$ which has played such a central role here in section 2; and partly to introduce some half-developed but intriguing possibilities. For better or worse, we organise the survey as follows.

Section 4.1 takes stock of where we stand with formal symbols and formal function elements. In section 4.1.1 we focus on formal symbols and isolate a number of mysterious aspects of their behaviour which we need to respect when working with them - the need to expand power series in both directions; the need to delay calculation involving symbols and retain stand-alone zeros which may be brought back to life before final evaluation; and the need to embed terms within a continuous "fabric" and use limits to resolve ambiguities as required. Section 4.1.2 extends the discussion to formal function elements and considers an extra subtlety arising from the fact that resolution into functions on evaluation can introduce delta-functions and other distributions. In section 4.1.3 we make the briefest of asides to note the connection of these ideas to Umbral calculus.

Section 4.2 extends our work in section 2 using formal quantities to compactly recast existing results and extend them. In section 4.2.1 we derive a simple formal expression for the Hurewicz zeta function, $\zeta_{H}(z_{0};s)$ and we use it to read off a useful identity for the formal symbol $\tau$ as well as to deduce power series behaviour of $\zeta_{H}$ both for $z_{0}$ near $0$ and as $z_{0}\rightarrow\infty$. In section 4.2.2 we then use this formal symbolic expression for $\zeta_{H}$ to demonstrate its close connection to the family of functions $\overset{\lor}{q}_{\rho}(\alpha)$, $\rho\in\mathbb{C}$.\footnote{This encompasses and extends the well-known connection between $\zeta_{H}$ and the Bernoulli polynomials.} Finally, in section 4.2.3 we show that the key Euler-McLaurin sum formula for the p-sum of the zeta function, on which we have relied so heavily already in these introductory papers, can also be expressed much more simply and compactly using formal symbols. This re-expression allows us to derive much more easily the critical result (equation \ref{psumzeta_X_q_cech_exact_2}) involving this p-sum function and the function element $\overset{\lor}{q}$ which was the key step in proving theorem 2. We conclude with brief comments drawing together all our different formal simplifications and noting some useful tools they now facilitate, such as change of "subject".

Section 4.3 then explores how our formal symbols and formal function elements behave under the operations of calculus - integration and differentiation. In some cases (section 4.3.1) this leads to surprisingly simple proofs using formal quantities of results which are otherwise highly non-trivial. In others (section 4.3.2) it vividly demonstrates in action many of the mysterious behavioural features of formal quantities enumerated in sections 4.1.1 and 4.1.2. Concluding comments draw together some of these characteristics, and also look ahead to future work\footnote{In the third of the sets of papers following on from this introductory set - the key idea introduced there is the generalised C\'{e}saro re-definition of the Mellin transform and resulting concept of a so-called "TLA-coefficient function"} which simultaneously connects power series behaviour of functions near $0$ and near $\infty$.

In section 4.4 we return to round out our analysis of the family of functions $\overset{\lor}{q}_{\rho}(\alpha)$ which has been so central in this paper. We show in section 4.4.1 how delta-functions or their derivatives arise in $\overset{\lor}{q}_{\rho}(\alpha)$ for each $\rho\in\mathbb{Z}_{<0}$ by considering how $\overset{\lor}{q}_{\rho}(\alpha)$ behaves as $\rho$ approaches such points and connecting to H\"{o}rmander's treatment of such distributions in [3]. In section 4.4.2 we consider the functions $\overset{\lor}{q}_{\rho}(\alpha)$ for arbitrary $\rho\in\mathbb{C}\setminus\mathbb{Z}_{<0}$ from a Fourier perspective as period-1 functions in $\alpha$. We calculate their Fourier coefficients and resultant Fourier series expression, and by recalling our formal-symbol formula for $\overset{\lor}{q}_{\rho}(\alpha)$ and specialising to $\alpha=0$ we derive the functional equation of $\zeta$ as an immediate corollary.

We conclude the paper in section 4.5 with a small miscellany of additional results and observations. These are of varying degrees of speculativeness but all of sufficient interest, in our opinion, to warrant being outlined - and we hope they will spark interest in further investigation. In section 4.5.1 we return to the Euler-McLaurin sum formula. Having used the formal symbol $\tau$ to re-express its application for the p-sum function of $\zeta$ in much more compact form (equation \ref{psumzeta_k_tau} in section 4.2.3), we show that a corresponding simplification can be obtained in general for the p-sum function of any amenable function $f$. The resulting formula re-expresses the general Euler-McLaurin sum formula in startlingly compact form. Sections 4.5.2 and 4.5.3 then present some additional stray observations regarding the behaviour of the $\overset{\lor}{q}_{\rho}(\alpha)$, $\rho\in\mathbb{C}$, and some surprising places they turn up. In section 4.5.2 we assess the ray-behaviour of the $\overset{\lor}{q}_{n}(\alpha)$, $n\in\mathbb{Z}_{\geq{0}}$ on the positive critical line in $\mathbb{C}$, and note yet another connection between these functions and the Gamma function on this critical line (this time relating to the coefficients in an asymptotic expansion as $T\rightarrow\infty$). We also show how the $\overset{\lor}{q}_{n}(\alpha)$ turn up in the asymptotic expansion of the binomial coefficient $\binom{\rho}{j}$. In section 4.5.3 we consider the functions $\overset{\lor}{q}_{\rho}(\alpha)$ spectrally under the application of the operator $\alpha \frac{\textrm{d}}{\textrm{d}\alpha}$ and its powers. This leads to interesting formulae giving $\left(\alpha \frac{\textrm{d}}{\textrm{d}\alpha}\right)^{\nu}$ in terms of a formal operator, $\epsilon$, and a modified form of Pascal's triangle adapted to these computations. This adds the notion of a \textit{formal operator} to our collection of formal objects and seems worth doing for this reason - even though the calculations are only partially developed and have not yet led to new results.

Oh I say old chap, I nearly forgot! There is also an appendix. It contains in section 5.1 the derivation of the asymptotic formula for $\binom{\rho}{j}$ as $j\rightarrow\infty$ which we promised in section 4.2.1. In section 5.2 we then conclude with another highly speculative half-calculation. This considers again the Fourier coefficients, $a_{n}(\rho)$, of $\overset{\lor}{q}_{\rho}(\alpha)$ calculated in section 4.4.2. By considering the integral defining $a_{n}(\rho)$ geometrically and in a generalised C\'{e}saro manner we derive an alternative formula for it. Although this computation is incomplete we think it is intriguing enough to include as  a parting shot in this appendix. This is because it shows that many of the key elements of the correct formula arise naturally in a generalised geometric C\'{e}saro context; and by comparison with our existing formula it implies a limiting integral result for the $\overset{\lor}{q}_{\rho}(\alpha)$ which seems interesting. What ho, rather!

Note that, broadly speaking, the subsections of section 4 are independent - and the reader may peruse some, all (or none) of them according to taste. The organization adopted above may represent a poor approach, but we are not sure what are the correct paths to follow. Our belief is simply that these formal jungles are both fascinating and useful. As such our hope is to encourage others to explore more deeply along these few existing pathways - and perhaps to uncover other paths waiting to be discovered - in the hope of finding a more coherent structure and set of tools hidden therein, maybe even a mathematical Machu Pichu beneath the canopy.

\section{A brief survey of the formal jungles}

\subsection{Formal symbols and formal function elements}

Formal symbols and formal function elements are clearly mysterious creatures. Let us count some of the ways.

\subsubsection{Formal symbols - mysteries and rules for calculation}

We have seen enough already in section 2\footnote{and will see much more in section 4.2} to conclude that formal symbols like $\tau$ are useful but require care in handling:\\
\\
\textbf{(i)} To begin with, expressions involving formal symbols should only be evaluated "at the end". Intermediate evaluation immediately leads to errors, for the simple reason that we need $\tau^{r}\tau^{s}=\tau^{r+s}$ but it is obviously \textit{not} the case that $\zeta(-r)\zeta(-s)=\zeta(-(r+s))$. Thus for example, in calculations using functions defined via formal symbols, such as $\overset{\lor}{q}_{\rho}(\alpha):=(\alpha-\tau)^{\rho}$ - where the coefficients of powers of $\alpha$ are expressions involving $\tau$ - it is critical to apply the binomial theorem, and combine fully algebraically, before undertaking any evaluation of powers of $\tau$. That is to say, we need to resolve all functions defined in terms of $\tau$ into powers of $\tau$ (or powers of $\tau$ times integer powers of $\ln(\tau)$) using Taylor series or other power series expansions, and then combine them all fully, before undertaking final evaluation.\\
\\
\textbf{(ii)} For $\overset{\lor}{q}_{n}(\alpha)$ we have seen that it is useful to consider not just $(\alpha-\tau)^{n}$ defined as $\sum_{j=0}^{n}(-1)^{j}\binom{n}{j}\alpha^{n-j}\tau^{j}$, but to also include the $j=-1$ term "to the left". This is what gives us the additional term  $\frac{\alpha^{n+1}}{n+1}$ and hence lets us write $\overset{\lor}{q}_{n}(\alpha)$ as $(\alpha-\tau)^{n}$. For this reason, in this paper we have in fact taken $(\alpha-\tau)^{n}$ to be defined as $\sum_{j=-1}^{n}(-1)^{j}\binom{n}{j}\alpha^{n-j}\tau^{j}$. Similarly for arbitrary power $\rho$, we take $(\alpha-\tau)^{\rho}:=\sum_{j=-1}^{\infty}(-1)^{j}\binom{\rho}{j}\alpha^{\rho-j}\tau^{j}$, and we likewise apply the same extension for other binomial expressions like $(z+\tau)^{\rho}$ and $(k+\tau)^{\rho}$ (both of which will turn up in section 4.2).

In fact, here we are really including not just the $j=-1$ term, but all the terms "to the left", so that $(\alpha-\tau)^{\rho}:=\sum_{j=-\infty}^{\infty}(-1)^{j}\binom{\rho}{j}\alpha^{\rho-j}\tau^{j}$. The additional terms for $j\leq{-2}$ are all what we call "stand-alone zero", i.e. as they stand, they all evaluate to zero since for $j=-2, -3, \ldots$ the zero in $\binom{\rho}{j}$ no longer has a countervailing pole in $\tau^{j}$ as it did when $j=-1$. More on these "stand-alone zeros" in point (iv) below.\\
\\
\textbf{(iii)} In making sense of this extra, non-zero term at $j=-1$ we need to consider the coefficient of $\alpha^{\rho+1}$ as the limit $\lim_{\epsilon\rightarrow 0}e^{i\pi (-1+\epsilon)}\binom{\rho}{-1+\epsilon}\tau^{-1+\epsilon}$. We are thus thinking of these discrete coefficients, indexed by $j\in\mathbb{Z}$, as actually embedded in a formal "coefficient fabric", with coefficient of $\alpha^{\rho-\nu}$ at arbitrary $\nu\in\mathbb{C}$ given by $e^{i\pi \nu}\binom{\rho}{\nu}\tau^{\nu}$.\\
\\
\textbf{(iv)} Overall we see that interesting functions (e.g. our rescaled versions of the Bernoulli polynomials $\overset{\lor}{q}_{n}(\alpha)$, their generalisations $\overset{\lor}{q}_{\rho}(\alpha)=(\alpha-\tau)^{\rho}$, and many more we will encounter in section 4.2) are given by well-defined power series whose coefficients involve formal symbols like $\tau$.

Since these must only be evaluated "at the end" we need to be very careful in future calculations not to omit terms which are zero when considered on a purely stand-alone basis, since they may generate non-zero contributions when they are combined with $\tau$-based coefficients arising elsewhere in the computation before evaluation. And the resultant contributions may then be determined by limiting calculations in which the symbol-based quantities at integer-indices are approached as limits through non-integer-indexed values in associated fabrics.\\
\\
\textbf{(v)} In a future paper, this fascinating business of retaining "stand-alone zeros" from power series extended in both directions until "the end", and approaching final evaluation via complex limits if required, will be shown to be both natural and useful. The theoretical rationale for it, however, is best deferred until our third set of additional papers following this introductory set. There we will demonstrate that the generalised C\'{e}saro framework is the natural framework for analysing Mellin transforms, and derive a fundamental result linking such Mellin transforms to a fabric which \textit{simultaneously} carries information about the two power series describing the underlying function near $0$ and near $\infty$.\\
\\
\textbf{(vi)} \textbf{[Summary of rules]:} For now, however, we ignore such theoretical questions and focus on the practicalities of computations with formal symbols. For this the take-away messages can be summarised as follows:

\textbf{(a)} Formal symbols can be combined with real or complex variables to express functions, with the symbols appearing in the coefficients of the power series for the function near $0$ or near $\infty$;

\textbf{(b)} In such cases, we must be disciplined not to evaluate early and instead wait "until the end" after all expressions in these symbols have been combined;

\textbf{(c)} These symbol-based coefficients in such power series lie within a whole coefficient "fabric" giving values for these coefficients not just at integer-indices, but at arbitrary complex values;

\textbf{(d)} The values of such coefficients at an integer-index must then sometimes be evaluated as a limit as we approach the integer-index through complex values in our coefficient fabric;

\textbf{(e)} Because of the need not to evaluate early and to consider fabric limits, we need to be careful not to ignore terms which are zero on a stand-alone basis ("retain the stand-alone zeros!"). As part of this, we may need to consider binomial and other power series expansions (e.g. Taylor series or asymptotic series near $0$ or $\infty$) as extending beyond the usual index-ranges which we are used to considering when all the variables involved are simply real or complex variables, to also encompass additional terms arising from index-values "to the left" or "to the right".

\subsubsection{Formal function elements}

If all of the above regarding formal symbols sounds a bit like voodoo - apologies! Alas, as regards formal function elements, there is even a further complication as follows:\\
\\
\textbf{(i)} In equation \ref{psumzeta_X_q_cech_exact_2} the expression $\left(X-\overset{\lor}{q}\right)^{-s}$ involving the formal function element $\overset{\lor}{q}$ arose naturally, but we interpreted it as $\sum_{j=0}^{\infty}(-1)^{j}\binom{-s}{j}X^{-s-j}\overset{\lor}{q}^{j}$, leaving the term $\frac{X^{1-s}}{1-s}$ as a separate term. What would happen if we instead try to include it as the $j=-1$ term in the above sum, in analogy with what we did for the interpretation of $(\alpha-\tau)^{\rho}$?

Well the $j=-1$ term gives the expression $-\binom{-s}{-1}\overset{\lor}{q}^{-1} \cdot X^{1-s}=-\binom{-s}{-1}\overset{\lor}{q}_{-1}(\alpha) \cdot X^{1-s}$. Now, to understand $\overset{\lor}{q}_{-1}(\alpha)$, recall that
\begin{equation}
\frac{d}{d\overset{\lor}{q}}\overset{\lor}{q}^{n}=n\overset{\lor}{q}^{n-1}=n\overset{\lor}{q}_{n-1}(\alpha)=\frac{d}{d\alpha}\overset{\lor}{q}_{n}(\alpha) \quad .
\label{q_cech_differentiation_1}\end{equation}
Since $\overset{\lor}{q}_{0}$ is the period-1 function given by $\overset{\lor}{q}_{0}(X)=\overset{\lor}{q}_{0}(\alpha)=\alpha-\frac{1}{2}$, it follows that formally we have
\begin{equation}
0\cdot\overset{\lor}{q}_{-1}(X)=1-\sum_{j=1}^{\infty}\delta_{j}(X)
\label{q_cech_minus_1}\end{equation}
where $\delta_{j}(X)$ is the delta-function located at $X=j$; and then
\begin{eqnarray}
(-1)\cdot 0\cdot\overset{\lor}{q}_{-2}(X) & = & -\sum_{j=1}^{\infty}\delta_{j}^{\prime}(X) \nonumber\\
\textrm{and} \quad (-2)\cdot (-1)\cdot 0\cdot\overset{\lor}{q}_{-3}(X) & = & -\sum_{j=1}^{\infty}\delta_{j}^{\prime \prime}(X) \quad \textrm{etc} \quad .
\label{q_cech_minus_2}\end{eqnarray}
If we seek to include the $j=-1$ term in our expression for $\left(X-\overset{\lor}{q}\right)^{-s}$ in equation \ref{psumzeta_X_q_cech_exact_2} and try to evaluate it by considering $(-1)^{-1+\epsilon}\binom{-s}{-1+\epsilon}\overset{\lor}{q}_{-1+\epsilon}(X) \cdot X^{-s+1-\epsilon}$ in the limit as $\epsilon\rightarrow0$, then we \textit{can} give meaning to the formal equation \ref{q_cech_minus_1} in a way which ends up generating the $\frac{X^{1-s}}{1-s}$ term. However, when we do so the additional delta-function terms in equation \ref{q_cech_minus_1} then contribute a further term consisting of $-\sum_{j=1}^{\infty}\frac{j^{1-s}}{1-s}$, which is formally $\frac{\zeta(s-1)}{s-1}$.

Moreover, the parallel no longer holds with the interpretation of $(\alpha-\tau)^{\rho}$. There the inclusion of the $j=-1$ term naturally gave the desired additional contribution, but the $j=-2,-3, \ldots$ terms were all still zero - so that the inclusion of the $j=-1$ term could in fact be viewed as the inclusion of all $j\in\mathbb{Z}_{<0}$, turning our sum defining $(\alpha-\tau)^{\rho}$ into $\sum_{j=-\infty}^{\infty}(-1)^{j}e^{i\pi \rho}\binom{\rho}{j}\tau^{\rho-j}\alpha^{j}$. In this case the terms arising from $j=-2,-3, \ldots$ are \textit{not} all zero. Rather, using the same limiting approach for interpretation, they seem to give an infinite set of additional terms, all equal to $-\sum_{j=1}^{\infty}\frac{j^{1-s}}{1-s}=\frac{\zeta(s-1)}{s-1}$ except alternating in sign. Thus if we were to try to expand the sum defining $\left(X-\overset{\lor}{q}\right)^{-s}$ to be $\sum_{j=-\infty}^{\infty}(-1)^{j}\binom{-s}{j}X^{-s-j}\overset{\lor}{q}^{j}$ we \textit{would} end up including the $\frac{X^{1-s}}{1-s}$ term that is currently included separately in equation \ref{psumzeta_X_q_cech_exact_2}; but we would also end up including an infinite set of extra terms which, by turns, add and subtract $\frac{\zeta(s-1)}{s-1}$. How to make sense of this is not clear to us at the moment - although it should be noted that we take comfort in the fact that these extraneous terms are at least all independent of $X$.

For this reason we leave equation \ref{psumzeta_X_q_cech_exact_2} in its current form and continue - at least for now - to interpret $\left(X-\overset{\lor}{q}\right)^{-s}$ just as $\sum_{j=0}^{\infty}(-1)^{j}\binom{-s}{j}X^{-s-j}\overset{\lor}{q}^{j}$. We believe, however, that the exploration just sketched out is very suggestive and certainly that it warrants further investigation.

For now, however, we simply note that it implies there may be additional mysteries regarding formal function elements, beyond even those which we have enumerated for formal symbols. This is owing to the fact that delta-functions and other distributions can arise (as limits of functions) upon final evaluation in the case of formal function elements\footnote{In section 4.4.2 we shall see how the functions $\overset{\lor}{q}_{\rho}(\alpha)$ do in fact lead to delta functions at the integer points $X=1,2,3, \ldots$ in the limit as $\rho\rightarrow -1$, and similarly for derivatives of delta-functions as $\rho\rightarrow -2$, $\rho\rightarrow -3$ etc. We shall do this by making contact with H\"{o}rmander's characterisation of such distributions in [3].}, whereas formal symbols can only ever evaluate (even in limiting cases) to complex values.

\subsubsection{Umbral calculus?}

The notions of formal symbols and formal function elements are clearly related to the theory of Umbral calculus. That being said, it is not clear (at least to an Umbral novice like the author) that they can be merely subsumed within that theory as it stands - either as regards our applications of them above, or as regards the key properties we have identified concerning their behaviour - such as bi-directional extension of series, resuscitation of stand-alone zeros, at-the-end calculation and embedding within continuous fabrics.

In any case, while Umbral calculus invokes linear functionals in an effort to render calculations like these rigorous, our focus at this stage is instead unapologetically on the \textit{techniques} and \textit{mechanics} of such calculations. Our aim for now is just to be able to calculate efficiently and reliably rather than become sidetracked into machinery and rigour.

\subsection{Using formal quantities to compactly recast existing results and extend them}

In section 2 we have seen many instances where the use of formal symbols or formal function elements leads to much cleaner expression of results and thus appears to be the "right way" to express them. In this section we give further examples using these quantities to extend results, develop new ones, and express them all in highly compact form.

\subsubsection{The Hurewicz zeta function and a useful symbolic identity} 
Consider $\zeta_{H}(z_{0};s)=R_{+}[\tilde{z}^{-s}](z_{0})$ for generic $s$. Clearly $\zeta_{H}(0;s)=\zeta(s)$ and $\frac{d}{dz_{0}}\zeta_{H}(z_{0};s)=-sR_{+}[\tilde{z}^{-s-1}](z_{0})$. Thus $\frac{d}{dz_{0}}\zeta_{H}(z_{0};s)\big|_{z_{0}=0}=-s\zeta(s+1)$; and $\frac{d^{2}}{dz_{0}^{2}}\zeta_{H}(z_{0};s)\big|_{z_{0}=0}=(-s)(-s-1)R_{+}[\tilde{z}^{-s-2}](0)=s(s+1)\zeta(s+2)$; and so on. Thus the Taylor series for $\zeta_{H}$ in $z_{0}$ around $0$ is
\begin{eqnarray}
\zeta_{H}(z_{0};s) & = & \zeta(s)-s\zeta(s+1)z_{0}+\frac{s(s+1)}{2!}\zeta(s+2)z_{0}^{2} - \ldots \nonumber\\
 & = & \sum_{j=0}^{\infty}\binom{-s}{j}\zeta(s+j)z_{0}^{j} \quad .
\label{zeta_H_TS_1}\end{eqnarray}
Since, as noted before (see Appendix 5.1), $\binom{-s}{j}\sim j^{s-1}$ as $j\rightarrow\infty$, this Taylor series is absolutely convergent with radius of convergence $1$. It follows that, at least on the disk $\big|z_{0}\big| <1$, $\zeta_{H}$ is given in terms of $\tau$ by the very simple formula
\begin{equation}
\zeta_{H}(z_{0};s) = \left(z_{0}+\tau\right)^{-s} = \sum_{j=0}^{\infty}\binom{-s}{j}\zeta(s+j)\,z_{0}^{j}
\label{zeta_H_TS_2}\end{equation}
where here we are now taking $\tau$ as the "subject" of the binomial, so that we get descending powers $\tau^{-s}$, $\tau^{-s-1} \ldots$ and non-negative integer powers of $z_{0}$, as befitting a Taylor series.

Note here that for $s\notin\mathbb{Z}$ the descending powers of $\tau$ never include $\tau^{-1}$ and nor do the "omitted" powers corresponding to $j\in\mathbb{Z}_{<0}$ in equation \ref{zeta_H_TS_2}. As such we do not have to deal with any singularities in the terms $\binom{-s}{j}\zeta(s+j)z_{0}^{j}$ for $j\geq{0}$ and these terms are all identically zero for $j<0$, so that their omission is uncontroversial.

For $s=-n\in\mathbb{Z}_{\leq{0}}$, the descending powers $\tau^{n}$, $\tau^{n-1} \ldots$ do contain $\tau^{-1}$ but the infinite series in equation \ref{zeta_H_TS_1} becomes a finite series, with the $\tau^{-1}$ singularity being cancelled by a zero in the factor $\binom{n}{n+1}$ to leave the familiar term $\frac{z_{0}^{n+1}}{n+1}$. As such, $\zeta_{H}(z_{0};s) = \left(z_{0}+\tau\right)^{-s}$ is well-defined and non-singular at such points (and the putative terms for $j<0$ are again all zero and uncontroversially omitted).\\
\\
\textbf{Connection of behaviour near $0$ and near $\infty$:} For $s=n\in\mathbb{Z}_{>0}$, however, things become interesting. When $n=1$ the constant $j=0$ term becomes singular reflecting the fact that $R_{+}[\tilde{z}^{-1}](z_{0})$ has a non-C\'{e}saro-amenable log-divergence at all $z_{0}$. When $n\in\mathbb{Z}_{>1}$, the terms for $j\geq{0}$ are all well-defined and continue to give the correct Taylor series for $\zeta_{H}(z_{0};n)$ around $0$. But now, if we take the expansion $\sum_{j=0}^{\infty}\binom{-s}{j}\tau^{-s-j}\,z_{0}^{j}$ and formally extend it "to the left" to include terms with $j\in\mathbb{Z}_{<0}$ we pick up negative powers of $z_{0}$ - starting with $z_{0}^{-(n-1)}$ and descending from there.

In fact these additional terms, taken together, correctly give the asymptotic series expansion for $\zeta_{H}(z_{0};n)$ as $z_{0}\rightarrow\infty$ which we would obtain by applying the Euler-McLaurin sum formula to the p-sum function $\sum_{l=1}^{k}\frac{1}{(z_{0}+l)^{n}}$. The coefficient expression $\binom{-s}{j}\tau^{-s-j}$ in the expansion of $(z_{0}+\tau)^{-s}$ is in fact what we will later (in the last of the three sets of papers following this introductory set) call a TLA-coefficient function. And what we see here is a first example of a general result which we will derive in those future papers, namely that such TLA-coefficient functions \textit{simultaneously} carry information about the coefficients in the power series expansions for the underlying function both near $0$ and near $\infty$. In those papers we will, however, only be dealing with complex variables - the fact that here this behaviour continues to hold even though one of the variables is a formal symbol is itself an interesting and useful thing.\\
\\
\textbf{An identity for $\tau$:} We conclude this section with a useful formal symbol identity which follows from equation \ref{zeta_H_TS_2}. Taking $s\in\mathbb{C}\setminus{\mathbb{Z}_{>0}}$, we note that for $Re(s)<0$ we have immediately from the definition of $\zeta_{H}$ that $\zeta_{H}(-1;s)=\zeta_{H}(0;s)$; and this clearly extends by analytic continuation (e.g. via strip-wise C\'{e}saro extension) to all $s\in\mathbb{C}\setminus{\mathbb{Z}_{>0}}$. It follows that
\begin{equation}
\left(\tau-1\right)^{-s} = \tau^{-s} \quad \textrm{for all} \quad s\in\mathbb{C}\setminus{\mathbb{Z}_{>0}} \; .
\label{tau_identity_1}\end{equation}

\subsubsection{The connection between $\zeta_{H}$ and the $\overset{\lor}{q}_{\rho}$}

Recall our C\'{e}saro-adapted set of periodic functions $\{\overset{\lor}{q}_{\rho}(\alpha)\}$, which represents a generalisation of the Bernoulli polynomials from positive integer index to arbitrary complex $\rho$. Comparing equation \ref{zeta_H_TS_2} with our defining equation for these $\overset{\lor}{q}_{\rho}(\alpha)$ in equation \ref{Function_elements_def_1_rho}, which also converged absolutely for $\big|\alpha\big|<1$ when $\rho\notin\mathbb{Z}_{<0}$, we see that $\overset{\lor}{q}_{\rho}$ is in fact really just $\zeta_{H}$ in disguise, at least on the unit disk in $\mathbb{C}$:
\begin{equation}
\overset{\lor}{q}_{\rho}(\alpha) = e^{i\pi \rho}\zeta_{H}(-\alpha;-\rho) \quad \textrm{on} \quad \big|\alpha\big|<1 \quad \textrm{for any} \quad \rho\in\mathbb{C}\setminus{\mathbb{Z}_{<0}} \; .
\label{q_cech_rho_zeta_H}\end{equation}
The fact that $(\tau-1)^{-s} = \tau^{-s}$ in equation \ref{tau_identity_1} is then what allows us to conclude that $\overset{\lor}{q}_{\rho}(1) = \overset{\lor}{q}_{\rho}(0)$ for all $\rho\in\mathbb{C}\setminus{\mathbb{Z}_{<0}}$ and hence to extend $\overset{\lor}{q}_{\rho}(\alpha)$ periodically in $\alpha$ on $\mathbb{R}$, with period $1$ (in contrast with $\zeta_{H}$, which of course is not periodic in $z_{0}$).

\subsubsection{Euler-McLaurin and C\'{e}saro analysis for $\zeta$}

Continuing in the same vein as sections 4.2.1 and 4.2.2 we now give two further instances where the use of formal symbols and formal function elements leads to dramatic simplification in proofs we have given, and also beautifully simplifies the formulation of key results:\\
\\
\textbf{(i)} Just as the Taylor series for $\zeta_{H}(z_{0};s)$ around $0$ just simplified so dramatically using $\tau$, so the Euler-McLaurin relationship for $\sum_{j=1}^{k}j^{-s}$, on which we relied so heavily for all our generalised C\'{e}saro analysis of $\zeta$ in [1] and [2], likewise collapses to a very compact form. This Euler-McLaurin sum formula for the p-sum of $\zeta$ is given in equation \ref{psumzeta_1} and, on recalling that the Bernoulli numbers satisfy $B_{n}=(-1)^{n-1} \cdot n \cdot \zeta(1-n)$, we can rewrite this simply as:
\begin{equation}
s_{\zeta,s}(k+\alpha)\,=\,\sum_{j=1}^{k}j^{-s}\,=\,\zeta(s)\,-\,(k+\tau)^{-s} \quad .
\label{psumzeta_k_tau}\end{equation}
Here we are making $k$ the "subject" in the expansion of the binomial $(k+\tau)^{-s}$ and this is extended to be defined as $\sum_{j=-1}^{\infty}\binom{-s}{j}k^{-s-j}\tau^{j}$ so that, as before, we absorb the $\frac{k^{1-s}}{1-s}$ term in equation \ref{psumzeta_1} into the binomial using the $j=-1$ term.\\
\\
\textbf{(ii)} In turn we can then combine this compact expression for the p-sum function with the equally compact expression for our $\overset{\lor}{q}_{\rho}(\alpha)$ as $(\alpha-\tau)^{\rho}$ to deduce almost trivially the key equation \ref{psumzeta_X_q_cech_exact_2} used in proving theorem 2 - an equation which previously required extended algebraic derivation.

Writing $X=k+\alpha$ as usual we can write $k+\tau$ as $X-(\alpha-\tau)$. Then in equation \ref{psumzeta_k_tau} we get that
\begin{eqnarray*}
\sum_{j=1}^{k}j^{-s} & = & \zeta(s)\,-\,(X-(\alpha-\tau))^{-s} = \zeta(s) - \sum_{j=-1}^{\infty}(-1)^{j}\binom{-s}{j}X^{-s-j}(\alpha-\tau)^{j}\\
 & = & \zeta(s) - \sum_{j=-1}^{\infty}(-1)^{j}\binom{-s}{j}X^{-s-j}\overset{\lor}{q}_{j}(\alpha)\\
 & = & \zeta(s) + \frac{X^{1-s}}{1-s} - \left(X-\overset{\lor}{q}\right)^{-s}
\end{eqnarray*}
which is equation \ref{psumzeta_X_q_cech_exact_2}. Here, as discussed in section 4.1.2, the expression $(X-\overset{\lor}{q})^{-s}$ is taken as just $\sum_{j=0}^{\infty}(-1)^{j}\binom{-s}{j}X^{-s-j}\overset{\lor}{q}^{j}$ and we have split the term $\frac{X^{1-s}}{1-s}$ out separately.

In this derivation we have glossed over any complications which arise regarding the appearance of delta-functions in $\binom{-s}{-1}\overset{\lor}{q}_{-1}(\alpha)$. It nonetheless serves to illustrate forcefully how use of the formal symbol $\tau$ and the formal function element $\overset{\lor}{q}$ allows us to derive this key equation in a way which seems almost unreasonably simple and direct.\\
\\
\textbf{Comment:} One further point is worth mentioning. Equation \ref{psumzeta_k_tau} can in fact be viewed as a formal special case of equation \ref{zeta_H_TS_2} with $z_{0}=k$, since $\sum_{j=1}^{k}j^{-s}=\zeta(s)-\zeta_{H}(k;s)$. We say "formal" here because the expression $\zeta_{H}(z_{0};s)=(z_{0}+\tau)^{-s}$ has only been validated for $\big|z_{0}\big|<1$ whereas here $k\in\mathbb{Z}_{\geq{1}}$. However, having validated equation \ref{zeta_H_TS_2} for $\big|z_{0}\big|<1$, it is easy to analytically extend also to $\big|z_{0}\big|\geq{1}$ by re-interpreting $(z_{0}+\tau)^{-s}$  to change the subject of the binomial from $\tau$ to $z_{0}$. Then $(z_{0}+\tau)^{-s}=\sum_{j=-1}^{\infty}\binom{-s}{j}\tau^{j}z_{0}^{-s-j}$, giving an asymptotic series in $z_{0}$ around $\infty$ (matching the asymptotic series given by the Euler-McLaurin sum formula) rather than a Taylor series in $z_{0}$ around $0$. This is another way of viewing the discussion after equation \ref{zeta_H_TS_2} in section 4.2.1.

This capacity - when dealing with series expansions with coefficients involving formal symbols or formal function elements - to change perspective on what is the "subject" in the expansion is a useful thing. It allows us to go from local series for small $z$ near $0$ to asymptotic series for large $z$ near $\infty$. At all times, however, it needs to be undertaken with appropriate care in extending the index-ranges of such series either "to the left" or "to the right", beyond what we may be accustomed to (much more on this, as promised, in future papers).

\subsection{Calculus with formal symbols and formal function elements}

Beyond power series algebra, can we do other operations when dealing with expressions involving formal quantities? In particular, can we do calculus - differentiation and integration - on such expressions, whether with respect to other complex variables involved in these expressions, or on the formal quantities involved in them directly? And, if so, can we use calculus methods to further simplify existing theory and then extend it? This section shows that the answer to these questions is, in essence, "Yes" - and in the course of computations demonstrating this, we will encounter many of the subtleties we have mentioned regarding the manipulation and interpretation of such formal quantities.

\subsubsection{First integration and differentiation results for expressions with formal symbols and formal function elements}

We give one example each of a result involving integration and a result involving differentiation of expressions containing formal symbols.\\
\\
\textbf{(1) [Integration result]:} In [2] we showed that
\begin{equation*}
\int_{-1}^{0}\zeta_{H}(z_{0};s)\,dz_{0} \,=\, 0 \quad \textrm{for} \quad s\in\mathbb{C}\setminus\mathbb{Z}_{>0}
\end{equation*}
using an argument based on right-Riemann sums and the dilation-invariance of C\'{e}saro convergence. Using equation \ref{zeta_H_TS_2} expressing $\zeta_{H}$ in terms of $\tau$, this instead becomes a trivial result. We have
\begin{eqnarray*}
\int_{-1}^{0}\zeta_{H}(z_{0};s)\,dz_{0} & = & \int_{-1}^{0}(z_{0}+\tau)^{-s}\,dz_{0} = \left[\frac{(z_{0}+\tau)^{1-s}}{1-s}\right]_{-1}^{0}\\
 & = & \frac{1}{1-s}\{\tau^{1-s}-(\tau-1)^{1-s}\} = 0
\end{eqnarray*}
by equation \ref{tau_identity_1}.\\
\\
\textbf{(2) [Differentiation result]:} Going back to theorems 1 and 2, consider the expression $\left(X-\overset{\lor}{q}\right)^{n}=\sum_{j=0}^{n}(-1)^{j}\binom{n}{j}X^{n-j}\overset{\lor}{q}^{j}$ which arose in our proofs of them and which involves the formal function element $\overset{\lor}{q}$. Since the $\left\{\overset{\lor}{q}_{j}\right\}_{j=0}^{\infty}$ form a C\'{e}saro-adapted scale, we know that
\begin{equation*}
\frac{d}{dX}\overset{\lor}{q}^{j}\,=\,\frac{d}{dX}\overset{\lor}{q}_{j}(X)\,=\,j\overset{\lor}{q}_{j-1}(X)\,=\,\frac{d}{d\overset{\lor}{q}}\overset{\lor}{q}^{j} \quad \textrm{for} \quad j\in\mathbb{Z}_{\geq{1}}
\end{equation*}
while
\begin{equation*}
\frac{d}{dX}\overset{\lor}{q}^{0}\,=\,\frac{d}{dX}\overset{\lor}{q}_{0}(X)\,=\,1-\sum_{l=1}^{\infty}\delta_{l}(X) \; .
\end{equation*}
It follows that
\begin{eqnarray*}
\frac{d}{dX}\left(X-\overset{\lor}{q}\right)^{n} & = & \{X^{n}\left[1-\sum_{l=1}^{\infty}\delta_{l}(X)\right]+nX^{n-1}\overset{\lor}{q}^{0}\}\\
 & & +\sum_{j=1}^{n}(-1)^{j}\binom{n}{j}\{X^{n-j}j\overset{\lor}{q}^{j-1}+(n-j)X^{n-j-1}\overset{\lor}{q}^{j}\}
\end{eqnarray*}
and since, for any $0\leq{j}\leq{n-1}$ we have that the coefficient of $X^{n-j-1}\overset{\lor}{q}^{j}$ is
\begin{equation*}
(-1)^{j}\binom{n}{j}(n-j)+(-1)^{j+1}\binom{n}{j+1}(j+1) = (-1)^{j}\frac{n!}{j!}\left\{\begin{array}{cc}\frac{1}{(n-j-1)!}\\
 \\
 -\frac{1}{(n-j-1)!}\end{array}\right\} = 0
\end{equation*}
so we get
\begin{equation}
\frac{d}{dX}\left(X-\overset{\lor}{q}\right)^{n}=X^{n}\cdot\left[1-\sum_{l=1}^{\infty}\delta_{l}(X)\right] \quad .
\label{X_minus_q_cech_n_derivative_1}\end{equation}
In particular, away from the integer points on $[0,\infty)$ we have that
\begin{equation}
\frac{d}{dX}\left(X-\overset{\lor}{q}\right)^{n}=X^{n} \quad .
\label{X_minus_q_cech_n_derivative_2}\end{equation}
\\
\textbf{Comments: (i)} This is a striking result given how complicated the function $\left(X-\overset{\lor}{q}\right)^{n}$ is and the fact that $R(X):=\left(X-\overset{\lor}{q}\right)^{n}$ is a quasi-oscillatory function which is \textit{strongly} convergent to $0$ in a generalised C\'{e}saro sense under $P^{n+1}$. It shows that this oscillatory character arises only because of drops by $k^{n}$ at each integer point $X=k$ in $[0,\infty)$, with simple, monotonic behaviour between drops. In this regard it is reminiscent of the famous argument of the Riemann zeta function, $S(T)$, which exhibits the same sort of pattern between and across non-trivial roots and which also gives rise to a C\'{e}saro-adapted scale (as we will show in our 2nd set of papers following on from this introductory set).\\
\\
\textbf{(ii)} Note that, while surprising, equations \ref{X_minus_q_cech_n_derivative_1} and \ref{X_minus_q_cech_n_derivative_2} fit with equation \ref{psumzeta_X_q_cech_exact_2} when $s=-n$. In equation \ref{psumzeta_X_q_cech_exact_2} the p-sum function $\sum_{j=1}^{k}j^{n}$ is a step-function with derivative zero between integer points and a derivative  of $k^{n}\delta_{k}(X)$ at each integer point $X=k$ where a step occurs, so the LHS has derivative $X^{n}\cdot\sum_{l=1}^{\infty}\delta_{l}(X)$. But in light of equation \ref{X_minus_q_cech_n_derivative_1} this is also immediately what we get from differentiating the RHS w.r.t. $X$ in equation \ref{psumzeta_X_q_cech_exact_2}. Note also that equations \ref{X_minus_q_cech_n_derivative_1} and \ref{X_minus_q_cech_n_derivative_2} do generalize from $s=-n\in\mathbb{Z}_{\leq{0}}$ to arbitrary $s\in\mathbb{C}\setminus\mathbb{Z}_{>0}$, but we have focussed on the integer case in order to avoid infinite sums and provide a simpler illustration.

\subsubsection{Further differentiation of expressions with formal symbols and formal function elements - mysteries in action}

In the previous subsection we have seen examples showing that we can perform calculus on expressions involving formal symbols and formal function elements and that this can lead to new observations or greatly simplify derivations of existing results. In this subsection we consider a further setting for such calculus with expressions involving formal symbols and formal function elements, but one which leads in interesting new directions.

We will only focus on two examples and only show a cursory level of working, but this suffices to illustrate some new aspects of computation with formal quantities. It will also illustrate the criticality of two of the mysterious aspects regarding formal computations canvassed in section 4.1 - namely the need to defer evaluation "until the end" and the need to "retain stand-alone zeros" (in this case arising from extension of a traditional binomial expansion fully "to the left" beyond its normal index range).

Consider the family of functions $\overset{\lor}{q}_{\rho}(\alpha)$ where we know, for example, that $\overset{\lor}{q}_{0}(\alpha)=\alpha - \frac{1}{2}$ and $\overset{\lor}{q}_{1}(\alpha)=\frac{1}{2}\alpha^{2} - \frac{1}{2}\alpha+\frac{1}{12}$. Let us consider the behaviour of these functions $\overset{\lor}{q}_{\rho}(\alpha)$ as $\rho\rightarrow0$ and as $\rho\rightarrow1$, at least to first order, i.e. let us try to calculate $\frac{\partial}{\partial \rho}\overset{\lor}{q}_{\rho}(\alpha)\big|_{\rho=0}$ and $\frac{\partial}{\partial \rho}\overset{\lor}{q}_{\rho}(\alpha)\big|_{\rho=1}$, starting with the case of $\rho=0$.\\
\\
\textbf{Example calculation (1) [$\frac{\partial}{\partial \rho}\overset{\lor}{q}_{\rho}(\alpha)\big|_{\rho=0}$]:} Since $\overset{\lor}{q}_{\rho}(\alpha)=(\alpha-\tau)^{\rho}$, if we let $\rho=\epsilon$ small it follows that
\begin{eqnarray*}
\overset{\lor}{q}_{\epsilon}(\alpha) & = & e^{i\pi \epsilon}(\tau-\alpha)^{\epsilon} = \{1 + i\pi \epsilon + \ldots\}\{\tau^{\epsilon} - \binom{\epsilon}{1}\tau^{\epsilon-1}\alpha+\binom{\epsilon}{2}\tau^{\epsilon-2}\alpha^{2} - \ldots\} \\
 \\
 & = & \{1 + i\pi \epsilon + \ldots\}\left\{\begin{array}{cc}\left[\zeta(0)-\zeta^{\prime}(0)\epsilon + \ldots \right] - \epsilon\left[-\frac{1}{\epsilon}+\gamma + \ldots \right]\alpha\\
 \\
  + \frac{\epsilon(\epsilon-1)}{2!}\left[\zeta(2)-\zeta^{\prime}(2)\epsilon + \ldots \right]\alpha^{2} - \ldots \end{array} \right\} \\
 \\ 
 & = & \{1 + i\pi \epsilon + \ldots\}\left\{(\alpha-\frac{1}{2})+\left[\begin{array}{cc}-\zeta^{\prime}(0)-\gamma\alpha - \frac{\zeta(2)}{2}\alpha^{2}\\
 \\
  - \frac{\zeta(3)}{3}\alpha^{3} - \ldots \end{array}\right]\epsilon+ O(\epsilon^{2})\right\} \\
 \\
  & = & \overset{\lor}{q}_{0}(\alpha)+\left\{i\pi \overset{\lor}{q}_{0}(\alpha)+\left[-\zeta^{\prime}(0)-\gamma\alpha - \sum_{j=2}^{\infty}\frac{\zeta(j)}{j}\alpha^{j}\right]\right\}\epsilon + O(\epsilon^{2})
\end{eqnarray*}
so
\begin{equation}
\frac{\partial}{\partial \rho}\overset{\lor}{q}_{\rho}(\alpha)\big|_{\rho=0} = i\pi \overset{\lor}{q}_{0}(\alpha)+\left[-\zeta^{\prime}(0)-\gamma\alpha - \sum_{j=2}^{\infty}\frac{\zeta(j)}{j}\alpha^{j}\right] \quad .
\label{d_drho_q_Cech_rho_at_0_1}\end{equation}
On the other hand, differentiating w.r.t. $\rho$ gives us formally
\begin{eqnarray}
\frac{\partial}{\partial \rho}\overset{\lor}{q}_{\rho}(\alpha)\Big|_{\rho=0} & = & \frac{\partial}{\partial \rho}\left(e^{i\pi \rho}(\tau-\alpha)^{\rho}\right)\Big|_{\rho=0} \nonumber\\
 & = & \left\{i\pi e^{i\pi \rho}(\tau-\alpha)^{\rho} + e^{i\pi \rho}(\tau-\alpha)^{\rho}\ln(\tau-\alpha)\right\}\Big|_{\rho=0} \nonumber\\
 & = & i\pi \overset{\lor}{q}_{0}(\alpha) + (\tau-\alpha)^{0}\ln(\tau-\alpha) \; .
\label{d_drho_q_Cech_rho_at_0_2}\end{eqnarray}
Is equation \ref{d_drho_q_Cech_rho_at_0_2} actually equivalent to equation \ref{d_drho_q_Cech_rho_at_0_1}? Well, being careful not to evaluate early, let us start by trying to take $(\tau-\alpha)^{0}$ as just being given by the "usual" expansion only "to the right", namely
\begin{equation}
(\tau-\alpha)^{0} = \sum_{j=0}^{\infty}(-1)^{j}\binom{0}{j}\tau^{-j}\alpha^{j} = \binom{0}{0}\tau^{0}\alpha^{0}-\binom{0}{1}\tau^{-1}\alpha^{1}+\binom{0}{2}\tau^{-2}\alpha^{2}-\ldots
\label{binom_0_expansion_omit_standalone_zeros}\end{equation}
while we also take $\ln(\tau-\alpha)=\ln(\tau(1-\tau^{-1}\alpha))=\ln(\tau)+\ln(1-\tau^{-1}\alpha)$ as being given by the usual Taylor series to the right, namely
\begin{equation}
\ln(\tau-\alpha) = \ln(\tau) - \sum_{j=1}^{\infty}\frac{\tau^{-j}\alpha^{j}}{j} \; .
\label{ln_tau_minus_alpha_expansion}\end{equation}
Also, recalling the defining formal relation for $\tau$ - namely that $\tau^{s}:=\zeta(-s)$ - and differentiating w.r.t. $s$, let us interpret $\tau^{-s}\ln(\tau)$ as
\begin{equation}
\tau^{-s}\ln(\tau) = -\zeta^{\prime}(s) \; .
\label{tau_minus_s_ln_tau_interpretation}\end{equation}
Then, using this relation and combining terms fully from our series expansions in equations \ref{binom_0_expansion_omit_standalone_zeros} and \ref{ln_tau_minus_alpha_expansion} before evaluating in order to calculate $(\tau-\alpha)^{0}\ln(\tau-\alpha)$, we find that we get two contributions to the calculation of this product, each involving a singular term in $\epsilon$. The first is
\begin{eqnarray}
 &   & \sum_{j=0}^{\infty}(-1)^{j}\binom{0}{j}\left(\tau^{-j}\ln(\tau)\right)\alpha^{j} \nonumber\\
 & = & \left\{\binom{0}{0}\tau^{0}\alpha^{0}-\binom{0}{1}\tau^{-1}\alpha^{1}+\binom{0}{2}\tau^{-2}\alpha^{2}-\ldots\right\}\ln(\tau) \nonumber\\
 & = & -\zeta^{\prime}(0) + \lim_{\epsilon\rightarrow 0}\frac{1}{\epsilon}\,\alpha - \alpha + \alpha\ln(\alpha) + O(\epsilon)
\label{tau_minus_alpha_0_ln_tau_minus_alpha_contribution_1}\end{eqnarray}
and the second is
\begin{equation}
-\left\{\sum_{j=0}^{\infty}(-1)^{j}\binom{0}{j}\tau^{-j}\alpha^{j}\right\}\cdot\left\{\sum_{l=1}^{\infty}\frac{\tau^{-l}\alpha^{l}}{l}\right\} = \left\{\begin{array}{cc} - \lim_{\epsilon\rightarrow 0}\frac{1}{\epsilon}\,\alpha - \gamma\alpha\\
\\ - \alpha\ln \alpha - \sum_{j=2}^{\infty}\frac{\zeta(j)}{j}\alpha^{j} \end{array}\right\} \; .
\label{tau_minus_alpha_0_ln_tau_minus_alpha_contribution_2}\end{equation}
In carrying out these computations we have, of course, had to perform limit-based interpretations of expressions like $\binom{0}{1}\left(\tau^{-1}\ln(\tau)\right)\alpha^{1}=-\lim_{\epsilon\rightarrow0}\binom{0}{1+\epsilon}\zeta^{\prime}(1+\epsilon)\alpha^{1+\epsilon}$  and $\binom{0}{0}\tau^{-1}\alpha^{1}=\lim_{\epsilon\rightarrow0}\binom{0}{\epsilon}\zeta(1+\epsilon)\alpha^{1+\epsilon}$ - and it is these which have produced the singular terms involving $\lim_{\epsilon\rightarrow0}\frac{1}{\epsilon}\,\alpha$ in these contributions. Fortunately, the singular terms in these two contributions are offsetting, and so when we combine them in equation \ref{d_drho_q_Cech_rho_at_0_2} we get
\begin{equation}
\frac{\partial}{\partial \rho}\overset{\lor}{q}_{\rho}(\alpha)\big|_{\rho=0} = i\pi \overset{\lor}{q}_{0}(\alpha)+\left[-\zeta^{\prime}(0)-(\gamma+1)\alpha - \sum_{j=2}^{\infty}\frac{\zeta(j)}{j}\alpha^{j}\right]
\label{d_drho_q_Cech_rho_at_0_3}\end{equation}
Unfortunately, however, this is \textit{wrong}! It disagrees with the correct formula in equation \ref{d_drho_q_Cech_rho_at_0_1} owing to the presence of the extra anomalous $-\alpha$ term! Where have we made our mistake?

Well, it turns out we should \textit{not} have omitted the stand-alone zero terms "to the left" in the binomial expansion for $(\tau-\alpha)^{0}$ in equation \ref{binom_0_expansion_omit_standalone_zeros}. That is, we should also have included terms $\sum_{j=-1}^{-\infty}(-1)^{j}\binom{0}{j}\tau^{-j}\alpha^{j}$ and so taken
\begin{eqnarray}
(\tau-\alpha)^{0} & = & \sum_{j=-\infty}^{\infty}(-1)^{j}\binom{0}{j}\tau^{-j}\alpha^{j} \nonumber\\
 & = & \left\{\ldots + \binom{0}{-2}\tau^{2}\alpha^{-2} - \binom{0}{-1}\tau^{1}\alpha^{-1}\right\} \nonumber\\
 &   & + \left\{\binom{0}{0}\tau^{0}\alpha^{0}-\binom{0}{1}\tau^{-1}\alpha^{1}+\binom{0}{2}\tau^{-2}\alpha^{2}-\ldots\right\} \; .
\label{binom_0_expansion_with_standalone_zeros}\end{eqnarray}
When we do this, the extra contribution from $\{\sum_{j=-1}^{-\infty}(-1)^{j}\binom{0}{j}\tau^{-j}\alpha^{j}\}\ln(\tau)$ is still zero on final evaluation. However, the extra contribution arising from $\{\sum_{j=-1}^{-\infty}(-1)^{j}\binom{0}{j}\tau^{-j}\alpha^{j}\} \cdot \{-\sum_{l=1}^{\infty}\frac{\tau^{-l}\alpha^{l}}{l}\})$ leads to a countable subset of non-zero terms requiring evaluation using $\lim_{\epsilon\rightarrow0}$ in the usual way and providing, on final evaluation, an extra total contribution of
\begin{equation*}
\left\{\frac{1}{1\cdot2} + \frac{1}{2\cdot3} + \frac{1}{3\cdot4} + \ldots\right\}\alpha = \alpha \; .
\end{equation*}
Inclusion of this extra contribution cancels out the anomalous $-\alpha$ in equation \ref{d_drho_q_Cech_rho_at_0_3} and leaves us with the correct final version of this formula, namely
\begin{equation*}
\frac{\partial}{\partial \rho}\overset{\lor}{q}_{\rho}(\alpha)\big|_{\rho=0} = i\pi \overset{\lor}{q}_{0}(\alpha)-\zeta^{\prime}(0)-\gamma\alpha - \sum_{j=2}^{\infty}\frac{\zeta(j)}{j}\alpha^{j}
\end{equation*}
in agreement with equation \ref{d_drho_q_Cech_rho_at_0_1}.

Once again we thus see that we can successfully differentiate an expression involving a formal symbol. But this example demonstrates many of the interesting subtleties, mentioned previously, which need to be observed in doing so. We needed to not evaluate "until the end"; consequently we needed to retain "stand-alone zeros" until final evaluation since they came back to life when combined with terms involving $\tau$ from other expressions; and we needed to include terms from binomial or other series expansions outside our normal index-ranges, in this case further "to the left".\\
\\
\textbf{Example calculation (2) [$\frac{\partial}{\partial \rho}\overset{\lor}{q}_{\rho}(\alpha)\big|_{\rho=1}$]:} The corresponding calculation for $\rho=1$ confirms all these points and also confirms the correctness and sufficiency of including the extra terms corresponding to negative index-values in our binomial expansion of $(\tau-\alpha)^{1}$ (in this case, the terms $\sum_{j=-1}^{-\infty}(-1)^{j}\binom{1}{j}\tau^{1-j}\alpha^{j}\}$). Taking $\rho=1+\epsilon$ and expanding to $O(\epsilon)$ we get that
\begin{eqnarray*}
\overset{\lor}{q}_{1+\epsilon}(\alpha) & = & \overset{\lor}{q}_{1}(\alpha)+\left\{i\pi \overset{\lor}{q}_{1}(\alpha)+\left[\begin{array}{cc}\zeta^{\prime}(-1)-(\frac{1}{2}+\zeta^{\prime}(0))\alpha\\
\\ 
+ (\frac{1}{2}-\frac{\gamma}{2\cdot 1})\alpha^{2}\\ 
\\
- \sum_{j=2}^{\infty}\frac{\zeta(j)}{j(j+1)}\alpha^{j+1} \end{array}\right]\right\}\epsilon + O(\epsilon^{2})
\end{eqnarray*}
so that
\begin{equation}
\frac{\partial}{\partial \rho}\overset{\lor}{q}_{\rho}(\alpha)\big|_{\rho=1} \,=\, i\pi \overset{\lor}{q}_{1}(\alpha)+\left[\begin{array}{cc}\zeta^{\prime}(-1)-(\frac{1}{2}+\zeta^{\prime}(0))\alpha\\
\\ 
+ (\frac{1}{2}-\frac{\gamma}{2\cdot 1})\alpha^{2} - \sum_{j=2}^{\infty}\frac{\zeta(j)}{j(j+1)}\alpha^{j+1} \end{array}\right]
\label{d_drho_q_Cech_rho_at_1_1}\end{equation}
On the other hand
\begin{equation}
\frac{\partial}{\partial \rho}\overset{\lor}{q}_{\rho}(\alpha)\big|_{\rho=1} \,=\, \frac{\partial}{\partial \rho}\left(e^{i\pi \rho}(\tau-\alpha)^{\rho}\right)\big|_{\rho=1} = i\pi \overset{\lor}{q}_{1}(\alpha) - (\tau-\alpha)^{1}\ln(\tau-\alpha) \; .
\label{d_drho_q_Cech_rho_at_1_2}\end{equation}
If we use the expansion for $\ln(\tau-\alpha)$ given in equation \ref{ln_tau_minus_alpha_expansion} but include only the standard terms $\sum_{j=0}^{\infty}(-1)^{j}\binom{1}{j}\tau^{1-j}\alpha^{j}\}$ in the binomial expansion then, after performing the computations, invoking relationship \ref{tau_minus_s_ln_tau_interpretation} and cancelling offsetting $\frac{1}{\epsilon}$ limit-singularities, we get in equation \ref{d_drho_q_Cech_rho_at_1_2} that
\begin{equation*}
\frac{\partial}{\partial \rho}\overset{\lor}{q}_{\rho}(\alpha)\big|_{\rho=1} = i\pi \overset{\lor}{q}_{1}(\alpha)+\left[\begin{array}{cc}\zeta^{\prime}(-1)-(\frac{1}{2}+\zeta^{\prime}(0))\alpha\\
\\ 
+ (\frac{3}{4}-\frac{\gamma}{2\cdot 1})\alpha^{2} - \sum_{j=2}^{\infty}\frac{\zeta(j)}{j(j+1)}\alpha^{j+1} \end{array}\right] \quad .
\label{d_drho_q_Cech_rho_at_1_3}\end{equation*}
This is again correct \textit{except} for having an erroneous extra term of $\frac{1}{4}\alpha^{2}$. But if we include the additional stand-alone zero terms at negative index in the binomial expansion for $(\tau-\alpha)^{1}$, namely $\sum_{j=-1}^{-\infty}(-1)^{j}\binom{1}{j}\tau^{1-j}\alpha^{j}$, then in the same fashion as before we get an extra contribution of
\begin{equation*}
-\left\{\frac{1}{3\cdot2\cdot1} + \frac{1}{4\cdot3\cdot2} + \frac{1}{5\cdot4\cdot3} + \ldots\right\}\alpha^{2} = -\frac{1}{4}\alpha^{2} \; .
\end{equation*}
This removes the anomaly and brings our calculation back into agreement with equation \ref{d_drho_q_Cech_rho_at_1_1}.\\
\\
\textbf{Comments:} Two final points regarding the pair of computations in this subsection:\\
\\
\textbf{(i)} The reader will recognize much of the expression for $\frac{\partial}{\partial \rho}\overset{\lor}{q}_{\rho}(\alpha)\big|_{\rho=0}$ in equation \ref{d_drho_q_Cech_rho_at_0_1} from the Taylor series for $\ln(\Gamma(1-\alpha))$ around $0$. Equation \ref{d_drho_q_Cech_rho_at_0_1} can be re-expressed as
\begin{equation}
\frac{\partial}{\partial \rho}\overset{\lor}{q}_{\rho}(\alpha)\big|_{\rho=0} = i\pi \overset{\lor}{q}_{0}(\alpha)-\zeta^{\prime}(0)-\ln(\Gamma(1-\alpha))
\label{d_drho_q_Cech_rho_at_0_1_alt}\end{equation}
and then expression \ref{d_drho_q_Cech_rho_at_1_1} for $\frac{\partial}{\partial \rho}\overset{\lor}{q}_{\rho}(\alpha)\big|_{\rho=1}$ can be re-expressed in terms of an anti-derivative of this, together with some additional correction terms. These seem like interesting observations in their own right.\\
\\
\textbf{(ii)} The reader may also ask at this point why it is that we needed to extend the binomial expansion of $(\tau-\alpha)^{\rho}$ all the way "to the left" to include negative index-values, but did not attempt anything similar for the expansion of $\ln(1-\tau^{-1}\alpha)$, contenting ourselves with just taking $\sum_{j=1}^{\infty}$ in equation \ref{ln_tau_minus_alpha_expansion}?

Apart from the fact that this makes things work out correctly, whereas extending to the left in equation \ref{ln_tau_minus_alpha_expansion} seems to cause problems (e.g. the next term $\frac{1}{0}\tau^{0}\alpha^{0}$ would appear to combine with other finite terms to leave singularities with no countervailing zero factor ...), why is it that extension is correct in the one case and not in the other? What precisely are the rules here for where to begin and end such expansions in general?

As hinted at earlier, the key feature of at least a first answer to this question relates to whether such terms are actually "attached" to $0$ (i.e. arise from a Taylor series/asymptotic series around $0$) or "attached" to $\infty$ (i.e. arise from a Taylor series/asymptotic series around $\infty$). 

In a sense the terms from $j=-1,-2, \ldots$ in extending equation \ref{ln_tau_minus_alpha_expansion} for $\ln(1-\tau^{-1}\alpha)$ actually reflect behaviour for $\alpha$ near $\infty$ rather than $0$; whereas the terms for $j=-1,-2, \ldots$ in extending the binomial expansion of $(\tau-\alpha)^{\rho}$ still relate to behaviour near $0$, albeit that the coefficients of the powers $\frac{1}{\alpha}$, $\frac{1}{\alpha^{2}}$ etc are (thankfully!) all zero at least on a stand-alone basis.

However, going into any greater detail here would both take us too far afield for purposes of this paper, and is also, in any case, best deferred until the future set of papers we have mentioned now several times. There, as promised, we will show that the generalised C\'{e}saro framework is in fact the natural setting for analysis of the Mellin transform, and hence derive a key result interpreting such Mellin transforms in terms of the \textit{simultaneous} behaviour of Taylor series/asymptotic series for the underlying function at \textit{both} $0$ and $\infty$.

So - no more on this for now! But hopefully the reader finds this very truncated discussion intriguing and is inspired to have a gander at these future papers.

\subsection{The family of functions $\{\overset{\lor}{q}_{\rho}(\alpha)\}$}

In this section we return to the inhabited fringes of the formal jungles and consolidate our understanding of the family of functions $\{\overset{\lor}{q}_{\rho}(\alpha)\}_{\rho\in\mathbb{C}}$. We have seen that for $\rho=n\in\mathbb{Z}_{\geq{0}}$ these are re-scaled versions of the Bernoulli polynomials, and for $\rho\in\mathbb{C}\setminus\mathbb{Z}_{<0}$ the Taylor series $\overset{\lor}{q}_{\rho}(\alpha)\,=\,e^{i\pi \rho}(\tau-\alpha)^{\rho}\,=\,e^{i\pi \rho}\sum_{j=0}^{\infty}(-1)^{j}\binom{\rho}{j}\tau^{\rho-j}\alpha^{j}$ defines $\overset{\lor}{q}_{\rho}(\alpha)$ locally uniformly in $\rho$ as an absolutely convergent power series with radius of convergence $\big|\alpha\big|<1$. Since $\overset{\lor}{q}_{\rho}(1)=\overset{\lor}{q}_{\rho}(0)$ it can be extended  as a period-1 function in $\alpha\in\mathbb{R}$, and for $\rho=-n\in\mathbb{Z}_{<0}$ we formally get equations \ref{q_cech_minus_1} and \ref{q_cech_minus_2} involving delta-functions and their derivatives at the points $X=j\in\mathbb{Z}_{>0}$ when considering $\overset{\lor}{q}_{-1}(\alpha)$, $\overset{\lor}{q}_{-2}(\alpha)$ etc. Let us consider this latter behaviour.

\subsubsection{How the delta-function behaviour of $\overset{\lor}{q}_{\rho}(\alpha)$ arises as $\rho\rightarrow -1$}

Recall from [3] the definition of the functions
\begin{equation}
\chi_{+}^{\rho}:=\frac{x_{+}^{\rho}}{\Gamma(\rho+1)} \quad \textrm{where} \quad x_{+}(x)=\begin{cases}
x & \quad \textrm{if} \; x\in\mathbb{R}_{\geq{0}}\\
0 & \quad \textrm{if} \; x\in\mathbb{R}_{<0}\end{cases}
\label{Hormander_chi_plus_defn}\end{equation}
H\"{o}rmander observes that $\frac{d}{dx}\chi_{+}^{\rho}=\chi_{+}^{\rho-1}$ and since $\chi_{+}^{0}$ is the one-sided Heaviside function, so $\chi_{+}^{-1}=\delta_{0}(x)$ and $\chi_{+}^{-2}=\delta_{0}^{\prime}(x)$, $\chi_{+}^{-3}=\delta_{0}^{\prime\prime}(x)$ etc. The delta-function $\delta_{0}(x)$ can thus be thought of as being a limiting quotient of a function that tends to a one-sided version of $\frac{1}{x}$, divided by a quantity (in $\rho$ not $x$) that tends to a simple pole with residue $1$; and similarly for $\delta_{0}^{\prime}(x)$, $\delta_{0}^{\prime\prime}(x)$ etc with higher one-sided reciprocal powers of $x$ and varying residues for the pole on the denominator.

With this perspective we can explain the development of the delta-function behaviour in $\overset{\lor}{q}_{\rho}(\alpha)$ as $\rho\rightarrow -1$, as captured formally in equation \ref{q_cech_minus_1} - and the similar behaviour captured in equation \ref{q_cech_minus_2} as $\rho\rightarrow -2$, $\rho\rightarrow -3$ etc. Recall (see Appendix 5.1) that we have locally uniformly in $\rho$ that
\begin{equation}
\binom{\rho}{j} = (-1)^{j}\cdot C \cdot j^{-\rho-1} \cdot \{1+\frac{a_{1}}{j}+\frac{a_{2}}{j^{2}}+\ldots\}
\label{rho_choose_j_asymptotic_1}\end{equation}
while, as $j\rightarrow\infty$, we have $\tau^{\rho-j}=\zeta(j-\rho)=1+2^{\rho-j}+3^{\rho-j}+\ldots \sim 1$ (ignoring the residual exponential decay). If we just retain the leading order, we see from our defining series that, up to an overall constant factor (depending only on $\rho$), we have
\begin{equation}
\overset{\lor}{q}_{\rho}(\alpha) \sim \sum_{j=0}^{\infty}j^{-\rho-1}\alpha^{j} \quad .
\label{q_cech_rho_leading_order}\end{equation}
Thus $\overset{\lor}{q}_{-1}(\alpha) \sim \sum_{j=0}^{\infty}\alpha^{j} \sim \frac{1}{1-\alpha}$ and this smoothly tends to $1$ as $\alpha\rightarrow 0^{+}$, but has a one-sided "$\frac{1}{x}$-style singularity" as $\alpha\rightarrow 1^{-}$ (i.e. if $\alpha=1-\epsilon$ then $\overset{\lor}{q}_{-1}(\alpha) \sim \frac{1}{\epsilon}$ as $\epsilon\rightarrow 0^{+}$). Since $\overset{\lor}{q}_{-1}(\alpha)$, like all the $\overset{\lor}{q}_{\rho}(\alpha)$, is periodic with period $1$, this one-sided "$\frac{1}{x}$-style singularity" is repeated from the left at $X=2,3,4, \ldots$. Equation \ref{q_cech_minus_1} then becomes exactly what one would expect from the H\"{o}rmander perspective, with the factor of $0$ in the formal statement representing the limit of a denominator approaching a simple pole as $\rho\rightarrow -1$ and combining with all these one-sided $\frac{1}{x}$-style singularities to give rise to the term $-\sum_{j=1}^{\infty}\delta_{j}(X)$ on the RHS.

In the same fashion we see that $\overset{\lor}{q}_{-2}(\alpha) \sim \sum_{j=0}^{\infty}j\alpha^{j} \sim \frac{\alpha}{(1-\alpha)^{2}}$ and $\overset{\lor}{q}_{-3}(\alpha) \sim \sum_{j=0}^{\infty}j^{2}\alpha^{j} \sim \frac{\alpha(1+\alpha)}{(1-\alpha)^{3}}$ and so on\footnote{In fact the leading order behaviour of $\overset{\lor}{q}_{-n}(\alpha)$ is given by $Li_{-n+1}(\alpha)$ where $Li_{j}$ is the $j^{th}$ polylogarithm function.}. Each of these is smooth from the right as we approach $0$; and, as $\alpha=1-\epsilon$ approaches $1$ from the left, we have $\overset{\lor}{q}_{-2}(\alpha) \sim \frac{1!}{\epsilon^{2}}$, $\overset{\lor}{q}_{-3}(\alpha) \sim \frac{2!}{\epsilon^{3}}$, $\overset{\lor}{q}_{-4}(\alpha) \sim \frac{3!}{\epsilon^{4}}$ and so on. With the periodic replication of these one-sided singularities at $X=2,3,4, \ldots$ we again see the H\"{o}rmander perspective becoming evident, leading to the results of equation \ref{q_cech_minus_2} and producing $-\sum_{j=1}^{\infty}\delta_{j}^{\prime}(X)$ on the RHS in the equation for $\overset{\lor}{q}_{-2}(X)$, and $-\sum_{j=1}^{\infty}\delta_{j}^{\prime\prime}(X)$ on the RHS in the equation for $\overset{\lor}{q}_{-3}(X)$, and so forth.

\subsubsection{The Fourier series for $\overset{\lor}{q}_{\rho}(X)$ and a corollary}

Since $\overset{\lor}{q}_{\rho}(X)=\overset{\lor}{q}_{\rho}(\alpha)$ is periodic on $\mathbb{R}$ with period $1$, it is natural to ask what are its Fourier coefficients and its resulting expression as a Fourier series. These Fourier expressions are well-known for $\rho=n\in\mathbb{Z}_{\geq{0}}$ when $\overset{\lor}{q}_{n}(\alpha)$ is just a re-scaling of the periodised Bernoulli polynomial $\tilde{B}_{n+1}(\alpha)$ - but we wish to derive them for arbitrary $\rho\in\mathbb{C}$. In this subsection we do this by classical methods (and in appendix 5.2 we also attempt this via a speculative alternative approach based on generalised geometric C\'{e}saro methods).

Let the Fourier series coefficients of $\overset{\lor}{q}_{\rho}(\alpha)$ be $\{a_{n}(\rho)\}_{n=-\infty}^{\infty}$. These are given by
\begin{equation}
a_{n}(\rho)=\int_{0}^{1}\overset{\lor}{q}_{\rho}(\alpha) e^{-2\pi in\alpha}\,d\alpha
\label{q_cech_FS_coefficients_defn_1}\end{equation}
with the resulting Fourier series reconstruction of $\overset{\lor}{q}_{\rho}$ being
\begin{equation}
\overset{\lor}{q}_{\rho}(\alpha)=\sum_{n=-\infty}^{\infty}a_{n}(\rho) e^{2\pi in\alpha} \quad .
\label{q_cech_FS_reconstruction_1}\end{equation}
Now from tables or elementary calculation we know that for $\overset{\lor}{q}_{0}(\alpha)=\alpha-\frac{1}{2}$ we have
\begin{equation*}
a_{n}(0)=\begin{cases}
\frac{1}{-2\pi in} & \quad \textrm{if} \; n\in\mathbb{Z}\setminus\{0\}\\
0 & \quad \textrm{if} \; n=0 \quad .\end{cases}
\end{equation*}
Since $\frac{d}{d\alpha}\overset{\lor}{q}_{m}(\alpha)=m\overset{\lor}{q}_{m-1}(\alpha)$, it follows by repeated integration by parts that for $\rho=m\in\mathbb{Z}_{>0}$
\begin{equation*}
a_{n}(m)=\begin{cases}
\frac{(-1)^{m}m!}{(-2\pi in)^{m+1}} & \quad \textrm{if} \; n\in\mathbb{Z}\setminus\{0\}\\
0 & \quad \textrm{if} \; n=0\end{cases}
\end{equation*}
and then in general that
\begin{equation}
a_{n}(\rho)=\begin{cases}
\frac{e^{i\pi \rho}\Gamma(\rho+1)}{(-2\pi in)^{\rho+1}} & \quad \textrm{if} \; n\in\mathbb{Z}\setminus\{0\}\\
0 & \quad \textrm{if} \; n=0 \quad .\end{cases}
\label{q_cech_FS_coefficients_1}\end{equation}
We thus have the following explicit Fourier reconstruction of $\overset{\lor}{q}_{\rho}(\alpha)$:
\begin{eqnarray}
\overset{\lor}{q}_{\rho}(\alpha) & = & \frac{e^{i\pi \rho}\Gamma(\rho+1)}{2^{\rho+1}\pi^{\rho+1}}\sum_{n=1}^{\infty}\left\{\frac{1}{(-i)^{\rho+1}n^{\rho+1}}e^{2\pi i n\alpha}+\frac{1}{(i)^{\rho+1}n^{\rho+1}}e^{-2\pi i n\alpha}\right\} \nonumber\\
 & = & \frac{e^{i\pi \rho}\Gamma(\rho+1)}{2^{\rho+1}\pi^{\rho+1}}\sum_{n=1}^{\infty}\left\{-\frac{e^{\frac{i\pi\rho}{2}}e^{2\pi in\alpha}}{i}+\frac{e^{-\frac{i\pi\rho}{2}}e^{-2\pi in\alpha}}{i}\right\}\frac{1}{n^{\rho+1}} \nonumber\\
 & = & -\frac{e^{i\pi \rho}\Gamma(\rho+1)}{2^{\rho}\pi^{\rho+1}}\sum_{n=1}^{\infty}\frac{\sin\left(2\pi n\alpha + \frac{\pi \rho}{2}\right)}{n^{\rho+1}} \qquad .
\label{q_cech_FS_reconstruction_2}\end{eqnarray}
As an immediate consequence of this note that since $\overset{\lor}{q}_{\rho}(\alpha)=e^{i\pi \rho}(\tau-\alpha)^{\rho}$ so $\overset{\lor}{q}_{\rho}(0)=e^{i\pi \rho}\tau^{\rho}=e^{i\pi \rho}\zeta(-\rho)$. It follows from equation \ref{q_cech_FS_reconstruction_2} that we have
\begin{equation*}
e^{i\pi \rho}\zeta(-\rho)=-\frac{e^{i\pi \rho}\Gamma(\rho+1)}{2^{\rho}\pi^{\rho+1}}\sin\left(\frac{\pi \rho}{2}\right)\zeta(\rho+1)
\end{equation*}
and on letting $\rho=-s$ and simplifying, this becomes precisely the functional equation of the Riemann zeta function, namely that
\begin{equation}
\zeta(s)=2^{s}\pi^{s-1}\sin\left(\frac{\pi s}{2}\right)\Gamma(1-s)\zeta(1-s) \quad .
\label{Riemann_zeta_functional_equation}\end{equation}
Thus the functional equation for $\zeta $ is simply one corollary of our general Fourier series expression for $\overset{\lor}{q}_{\rho}(\alpha)$ in equation \ref{q_cech_FS_reconstruction_2}, obtained by specialising to $\alpha=0$. A similar derivation occurs if we specialise to the point $\alpha=\frac{1}{2}$. This is not surprising since we have seen that $\overset{\lor}{q}_{\rho}(\frac{1}{2})$ is very closely related to $\zeta_{H}(-\frac{1}{2};-\rho)$ and we have shown in [2] how the functional equation for $\zeta$ arises from consideration of $\zeta_{H}$ at $z_{0}=-\frac{1}{2}$ - but we omit further discussion here.

\subsection{Miscellaneous further observations and results}

\subsubsection{A compact formal expression for the Euler-McLaurin sum formula in general}

We have seen in section 4.2.3 how the Euler-McLaurin sum formula for the p-sum function of the zeta function can be expressed in the compact form given in equation \ref{psumzeta_k_tau} using the formal symbol $\tau$. Is similar simplification achievable in the case of the Euler-McLaurin sum formula in general? The answer is "Yes".

For $f$ an arbitrary amenable function, the Euler-McLaurin sum formula says that
\begin{equation}
\sum_{n=1}^{k}\,f(n)\,\sim\,\int_{0}^{k}\,f(x)\,dx\,+\,C_{f}\,+\,\frac{1}{2}f(k)\,+\,\sum_{r=2}^{\infty}\,\frac{B_{r}}{r!}\,f^{\left(r-1\right)}(k) \quad .\label{EulerMcLaurin_1}\end{equation}
Now recalling that $B_{n}=(-1)^{n-1} \cdot n \cdot \zeta(1-n)$ for all $n\in\mathbb{Z}_{\geq{1}}$ we have that
\begin{equation*}
\frac{1}{2}f(k)\,+\,\sum_{r=2}^{\infty}\,\frac{B_{r}}{r!}\,f^{\left(r-1\right)}(k) = -\zeta(0)f(k) - \sum_{r=1}^{\infty}\,\frac{\zeta(-r)}{r!}\,f^{\left(r\right)}(k)
\end{equation*}
and this can be re-expressed formally using $\tau$ as
\begin{equation*}
\frac{1}{2}f(k)\,+\,\sum_{r=2}^{\infty}\,\frac{B_{r}}{r!}\,f^{\left(r-1\right)}(k) = - \sum_{r=0}^{\infty}\,\frac{\tau^{r}D^{r}}{r!}[f](k)
\end{equation*}
where $D=\frac{\textrm{d}}{\textrm{d}x}$ is the usual derivative operator. The expression $\sum_{r=0}^{\infty}\,\frac{\tau^{r}D^{r}}{r!}$looks like the usual expansion for $\textrm{e}^{\tau D}$, but if we also extend the Taylor series for the exponential to the left we pick up a single extra term in the usual way from the cancellation of the singularities in $\tau^{-1}$ and $(-1)!$ via limits. This extra term is $-D^{-1}[f](k)$ and represents the negative of an anti-derivative of $f$ evaluated at $k$, i.e. $-\int_{0}^{k}\,f(x)\,dx$.\footnote{Any difference by a constant arising in such an anti-derivative - for example by choosing a different lower limit on the integral - is simply absorbed into the constant term $C_{f}$} Thus, apart from the $C_{f}$-term, the whole of the expression on the RHS in equation \ref{EulerMcLaurin_1} can be subsumed formally as $-\textrm{e}^{\tau D}[f](k)$; and since $D$ is the generator of translations\footnote{i.e. $\textrm{e}^{h D}[f](x_{0})=f(x_{0}+h)$}, this can be further simplified as just $-f(k+\tau)$. Thus using the formal symbol $\tau$ we can re-express the general Euler-McLaurin sum formula in the following remarkably compact form:
\begin{equation}
\sum_{n=1}^{k}\,f(n)\,\sim\,C_{f}\,-\,f(k+\tau) \quad .
\label{EulerMcLaurin_2}\end{equation}\\
\textbf{Comment:} Equation \ref{EulerMcLaurin_2} is potentially well-adapted for applying C\'{e}saro methods since, exactly as in section 2, we can write $k+\tau=(X-\alpha)+\tau=X-(\alpha-\tau)$ and thereby use a suitable series expansion for $f$ to obtain a series expansion for $f(k+\tau)$ in terms of powers of $X$ and powers of $\overset{\lor}{q}$ (remembering that $\overset{\lor}{q}^{\rho}=\overset{\lor}{q}_{\rho}(\alpha)=(\alpha-\tau)^{\rho}$). As we saw in section 2, this may allow the application of theorems 1 and 2 and thus facilitate C\'{e}saro analysis of the p-sum.

\subsubsection{Ray behaviour of $\overset{\lor}{q}_{\rho}(\alpha)$ and further connection to the Gamma function}

A couple of observations regarding the behaviour of the $\overset{\lor}{q}_{\rho}(\alpha)$ and their relationship to other functions and expressions:\\
\\
\textbf{(a)} If we let $Cr_{pos}$ denote the positive "critical line" in $\alpha$ (not $\rho$), namely $Cr_{pos}=\{\alpha=\frac{1}{2}+it \mid t\geq{0}\}$, then the $\{\overset{\lor}{q}_{n}(\alpha)\}_{n\in\mathbb{Z}_{\geq{0}}}$ seem to exhibit "ray behaviour". By this we mean that they map $Cr_{pos}$ to a subset of the rays emanating from $0$ along the positive or negative real or imaginary axes, with the image-ray seeming to rotate smoothly as $n$ increases, i.e.
\begin{eqnarray*}
\overset{\lor}{q}_{0}\,:\,Cr_{pos} \longmapsto Im_{pos} \quad \textrm{with} \quad \overset{\lor}{q}_{0}(\frac{1}{2})=0\\
\overset{\lor}{q}_{1}\,:\,Cr_{pos} \longmapsto Re_{neg} \quad \textrm{with} \quad \overset{\lor}{q}_{1}(\frac{1}{2})=-\frac{1}{24}\\
\overset{\lor}{q}_{2}\,:\,Cr_{pos} \longmapsto Im_{neg} \quad \textrm{with} \quad \overset{\lor}{q}_{2}(\frac{1}{2})=0\\
\overset{\lor}{q}_{3}\,:\,Cr_{pos} \longmapsto Re_{pos} \quad \textrm{with} \quad \overset{\lor}{q}_{3}(\frac{1}{2})=\frac{7}{960}\\
\overset{\lor}{q}_{4}\,:\,Cr_{pos} \longmapsto Im_{pos} \quad \textrm{with} \quad \overset{\lor}{q}_{4}(\frac{1}{2})=0\\
\overset{\lor}{q}_{5}\,:\,Cr_{pos} \longmapsto Re_{neg} \quad \textrm{with} \quad \overset{\lor}{q}_{5}(\frac{1}{2})=-\frac{31}{8064}\\
\overset{\lor}{q}_{6}\,:\,Cr_{pos} \longmapsto Im_{neg} \quad \textrm{with} \quad \overset{\lor}{q}_{6}(\frac{1}{2})=0\\
\overset{\lor}{q}_{7}\,:\,Cr_{pos} \longmapsto Re_{pos} \quad \textrm{with} \quad \overset{\lor}{q}_{7}(\frac{1}{2})=\frac{(2^{7}-1)}{10!}\cdot\frac{945}{8}\\
\vdots
\end{eqnarray*}
This seems interesting behaviour and, in particular, it would be good to understand whether this continues to occur for $\rho\in\mathbb{R}$ between integer values, and indeed for $\rho\in\mathbb{C}$ more generally. The min and max values given by $\overset{\lor}{q}_{n}(\frac{1}{2})$ noted above are also interesting and are familiar. Adopting the notation of [4], if we define
\begin{equation}
\theta(T)\,:=\,Im\left(\ln\left(\Gamma\left(\frac{\frac{1}{2}+iT}{2}\right)\right)\right)-\frac{T}{2}\ln \pi
\label{Theta_T_definition}\end{equation}
then the counting function for the non-trivial zeros of $\zeta$, $N(T)$, is given by
\begin{equation}
N(T)\,=\,\frac{1}{\pi}\theta(T)+1+S(T)
\label{N_T_counting_function_1}\end{equation}
where $S(T)$ is the famous argument of the zeta function and $\theta(T)$ has the asymptotic expansion
\begin{equation}
\theta(T)\,=\,\frac{T}{2}\ln\left(\frac{T}{2\pi}\right)-\frac{T}{2}-\frac{\pi}{8}+\frac{1}{48T}+\frac{7}{5760T^{3}}+\frac{31}{80640T^{5}}+\ldots
\label{Theta_T_AS}\end{equation}
(see [5, sections 6.7 and 6.5]). We see that, up to sign, the values of $\overset{\lor}{q}_{n}(\frac{1}{2})$ for $n\geq{1}$ correspond to the coefficients of $\frac{1}{T^{n}}$ in this equation, with
\begin{equation}
\big| \overset{\lor}{q}_{n}\left(\frac{1}{2}\right) \big| \,=\, (2n)\cdot \left(\textrm{coefficient of} \; \frac{1}{T^{n}}\right) \; .
\label{q_cech-half_Theta_T_AS_coefficients}\end{equation}
This is another example (like equation \ref{d_drho_q_Cech_rho_at_0_1_alt}) showing close connection between the functions $\overset{\lor}{q}_{\rho}(\alpha)$ and the log of the $\Gamma$-function, at least at integer values of $\rho$.\\
\\
\textbf{(b)} The alternative polynomials given by $\overset{\lor}{q}_{n}(\alpha)+\tau^{n}=\overset{\lor}{q}_{n}(\alpha)+\zeta(-n)$ exhibit the same ray behaviour while now also having $\alpha=0$ and $\alpha=1$ as roots and retaining $\alpha=\frac{1}{2}$ as a root for $n\in\mathbb{Z}_{>0}$ even. As such, it may be preferable to instead work with this family of polynomials. For example, if we consider $\ln\binom{\rho}{j}$ then the $\overset{\lor}{q}_{n}(\alpha)+\tau^{n}$ turn up as coefficients in its formal asymptotic expansion as $j\rightarrow\infty$, namely\footnote{We omit derivation of this formula, but the working of appendix 5.1 contains the initial steps}
\begin{equation}
\ln \binom{\rho}{j} = \left\{\begin{array}{cc} i\pi j - (\rho+1)\ln j -\left\{\begin{array}{cc}\gamma(\rho+1)+\frac{1}{2}\tau^{-2}(\rho+1)^{2}\\
\\
+\frac{1}{3}\tau^{-3}(\rho+1)^{3} + \ldots \end{array}\right\}\\
\\
+\left\{\begin{array}{cc}\left(\overset{\lor}{q}_{1}(\rho+1)+\tau^{1}\right)\frac{1}{j}+\frac{1}{2}\left(\overset{\lor}{q}_{2}(\rho+1)+\tau^{2}\right)\frac{1}{j^{2}}\\
\\
+\frac{1}{3}\left(\overset{\lor}{q}_{3}(\rho+1)+\tau^{3}\right)\frac{1}{j^{3}}+\ldots \end{array}\right\} \end{array}\right\}
\label{ln_rho_choose_j_AS}\end{equation}
This expansion again seems interesting in and of itself, and is susceptible to further simplification using a formal function element corresponding to this alternative family of polynomials.

\subsubsection{A few more stray thoughts re the functions $\overset{\lor}{q}_{\rho}(\alpha)$}

In subsection 4.4.1 we explained the appearance of delta-function behaviour in the $\overset{\lor}{q}_{\rho}(\alpha)$ as $\rho\rightarrow n\in\mathbb{Z}_{<0}$ by considering $(\alpha-\tau)^{\rho}$ and looking at leading order to focus on expressions of the form $\sum_{j}j^{n}\alpha^{j}$. It is independently interesting to consider expressions of the form $\sum_{j}j^{\nu}\alpha^{j}$ (which are the polylogarithm functions) for $\nu\in\mathbb{C}$ arbitrary. Now, by spectral definition, each term in the sum is given by $j^{\nu}\alpha^{j}=\left(\alpha\frac{d}{d\alpha}\right)^{\nu}[\alpha^{j}]$ (since $\alpha^{j}$ is an eigenfunction of the operator $\alpha\frac{d}{d\alpha}$ with eigenvalue $j$), and this leads us to want to better understand the operator $\left(\alpha\frac{d}{d\alpha}\right)^{\nu}$. We omit all details here, but this leads in further interesting formal directions.

To begin with, it turns out we can express $\left(\alpha\frac{d}{d\alpha}\right)^{\nu}$ as a series in the formal operator $\epsilon$ defined by $\epsilon^{m}:=\alpha^{m}\left(\frac{d}{d\alpha}\right)^{m}$. Specifically, we have $\left(\alpha\frac{d}{d\alpha}\right)^{\nu}=\epsilon e^{(T-1)\epsilon}\left[z^{\nu}\right] \big|_{z=1}$ where $T$ is the translation operator in $z$ given by $T[g](z):=g(z+1)$ - and this allows us to express $\left(\alpha\frac{d}{d\alpha}\right)^{\nu}$ even for $\nu\notin\mathbb{Z}$ in terms solely of integer powers of $\epsilon$. So now we see that interesting simplification can occur not just by using formal symbols and formal function elements, but also by introducing formal \textit{operators}. 

When $\nu=n$ is an integer this formula of course reduces to a finite sum\footnote{In the formula this occurs because $(T-1)$ is the discrete derivative and acts on a polynomial in a way which reduces its degree by $1$, so that all sufficiently high powers of $(T-1)$ reduce $z^{n}$ to $0$ and the series expansion of $e^{(T-1)\epsilon}$ terminates after a finite number of terms in its action on $z^{n}$} and we get an interesting modified analogue of Pascal's triangle describing the coefficients of the powers of $\epsilon$:\footnote{here for example, the 5th row below says that $\left(\alpha\frac{d}{d\alpha}\right)^{5}[f]=\alpha^{5}\frac{d^{5}f}{d\alpha^{5}} + 10\alpha^{4}\frac{d^{4}f}{d\alpha^{4}} + 25\alpha^{3}\frac{d^{3}f}{d\alpha^{3}} + 15\alpha^{2}\frac{d^{2}f}{d\alpha^{2}} + \alpha\frac{df}{d\alpha}$ and so on}
\begin{equation}
\begin{array}{cc}
1 \\
1 \qquad 1 \\
1 \qquad 3 \qquad 1 \\
1 \qquad 6 \qquad 7 \qquad 1 \\
1 \qquad 10 \qquad 25 \qquad 15 \qquad 1 \\
1 \qquad 15 \qquad 65 \qquad 90 \qquad 31 \qquad 1 \\
1 \qquad 21 \qquad 140 \qquad 350 \qquad 301 \qquad 63 \qquad 1 \\
\vdots
\end{array}
\end{equation}   
Just as the binomial coefficients appearing in the usual Pascal's triangle are the sum of the two terms appearing immediately above them, the coefficients in this amended version are also given by a simple sum in terms of the two terms immediately above them. The only difference is that now in the sum, the term above and to the left is multiplied by a constant factor corresponding to the power of $\epsilon$ it is a coefficient of, before being added to the term above and to the right. 

Thus for example $7=2 \cdot 3+1$ and $15=2 \cdot 7+1$ and $31=2 \cdot 15 +1$ etc; while $25=3 \cdot 6+7$ and $90=3 \cdot 25+15$ and $301=3 \cdot 90+31$ etc; and $65=4 \cdot 10+25$ and $350=4 \cdot 65+90$ etc; and so on. 

Just as with the traditional Pascal's triangle, this relationship allows it to be built extremely efficiently without having to either undertake long-hand derivation of $\left(\alpha\frac{d}{d\alpha}\right)^{n}[f]$ or long-hand calculation from the formula above that $\left(\alpha\frac{d}{d\alpha}\right)^{n}=\epsilon e^{(T-1)\epsilon}\left[z^{n}\right] \big|_{z=1}$. And it also gives a way of thinking of the coefficients which occur in the infinite series when $\nu\notin\mathbb{Z}$ as being embedded in a continuous fabric underlying this discrete triangle (in the same way the traditional binomial coefficients extend from $\binom{n}{j}$ appearing in Pascal's triangle to general $\binom{\rho}{j}$ for $\rho\in\mathbb{C}$).

Whether the reader's interest is more analytic - relating to operators and the asymptotic and distributional behaviour of the $\overset{\lor}{q}_{\rho}(\alpha)$; or purely combinatorial - relating to this amended Pascal's triangle and formal expressions - we believe that these stray observations are worth pursuing further.

\section{Appendix}

\subsection{An asymptotic expansion for $\binom{\rho}{j}$ as $j\rightarrow\infty$}

In several places we have used the following result:\\
\\
\textbf{Result:} \textit{Fix} $\rho\in\mathbb{C}$ \textit{arbitrary. Then as} $j\rightarrow\infty$ \textit{we have the asymptotic expansion}
\begin{equation}
\binom{\rho}{j} = (-1)^{j} \cdot C \cdot j^{-\rho-1} \cdot \left\{1+\frac{a_{1}}{j}+\frac{a_{2}}{j^{2}}+\ldots\right\}
\end{equation}
\textit{for some constants} $C$, $a_{1}$, $a_{2} \ldots$. \textit{This expansion is locally uniform in a neighbourhood of} $\rho$.\\
\\
\textbf{Proof:} We have
\begin{eqnarray*}
\binom{\rho}{j} & = & \frac{\rho \cdot (\rho-1) \cdots (\rho-j+1)}{1 \cdot 2 \cdots j} \\
 \\
 & = & (-1)^{j} \cdot \left(\frac{j-(\rho+1)}{j}\right) \cdot \left(\frac{(j-1)-(\rho+1)}{(j-1)}\right) \cdots \left(\frac{1-(\rho+1)}{1}\right) \\
 \\
 & = & (-1)^{j} \cdot \left(1-\frac{(\rho+1)}{j}\right) \cdot \left(1-\frac{(\rho+1)}{(j-1)}\right) \cdots \left(1-\frac{(\rho+1)}{1}\right) \quad .
\end{eqnarray*}
Therefore
\begin{equation*}
\ln \binom{\rho}{j} = \pm \, i\pi j + \sum_{i=1}^{j}\ln\left(1-\frac{(\rho+1)}{i}\right) \quad .
\end{equation*}
To understand this asymptotically we apply the Euler-McLaurin sum formula to $f(x)=\ln\left(1-\frac{(\rho+1)}{x}\right)$. Note first that
\begin{equation*}
f^{\prime}(x)\big|_{x=j} = \frac{(\rho+1)}{j(j-(\rho+1))} = \frac{(\rho+1)}{j^{2}} \cdot \left\{1+\frac{c_{1}}{j}+\frac{c_{2}}{j^{2}}+\ldots\right\}
\end{equation*}
and similarly for the higher derivatives of $f$. On the other hand
\begin{eqnarray*}
\int^{j}\ln\left(1-\frac{(\rho+1)}{x}\right)\,\textrm{d}x & = & \left[x \ln\left(1-\frac{(\rho+1)}{x}\right)\right]^{j} - \int^{j}\frac{(\rho+1)}{x-(\rho+1)}\,\textrm{d}x \\
 \\
 & = & j \ln\left(1-\frac{(\rho+1)}{j}\right) - (\rho+1)\ln\left(j-(\rho+1)\right) \\
 \\
 & = & -(\rho+1)\ln j - (\rho+1) + \frac{d_{1}}{j}+\frac{d_{2}}{j^{2}}+\ldots
\end{eqnarray*}
Therefore, in the Euler-McLaurin sum formula
\begin{equation*}
\ln \binom{\rho}{j} = -(\rho+1)\ln j \pm \, i\pi j - (\rho+1) + \frac{b_{1}}{j}+\frac{b_{2}}{j^{2}}+\ldots
\end{equation*}
and thus
\begin{equation*}
\binom{\rho}{j} = (-1)^{j} \cdot C \cdot j^{-(\rho+1)} \cdot \left\{1 + \frac{a_{1}}{j}+\frac{a_{2}}{j^{2}}+\ldots\right\}
\end{equation*}
and the constant $C$ can be chosen locally uniformly to apply within a neighbourhood of $\rho$. This completes the proof.

\subsection{A speculative generalised geometric C\'{e}saro approach to calculating the Fourier coefficients of $\overset{\lor}{q}_{\rho}(\alpha)$?}

Just as we ended section 4 with some rather speculative stray observations and calculations, we end this appendix with another speculative attempted calculation. Although we have not fully succeeded with it, we believe it shows enough to suggest that there must be something "there" which is worth exploring further. Specifically, let us try to calculate the Fourier coefficients given in equation \ref{q_cech_FS_coefficients_1} using a different heuristic approach based on geometric generalised C\'{e}saro methods.

We start by considering the integral for $a_{n}(\rho)$ in equation \ref{q_cech_FS_coefficients_defn_1} as defining a value attached to the point $\rho\in\mathbb{C}$. Since both $\overset{\lor}{q}_{\rho}(\alpha)$ and $\textrm{e}^{-2\pi in\alpha}$ are periodic with period $1$, it follows by integration by parts that we have
\begin{equation*}
a_{n}(\rho) = \left(\frac{2\pi in}{\rho+1}\right) \cdot \int_{0}^{1}\overset{\lor}{q}_{\rho+1}(\alpha) e^{-2\pi in\alpha}\,d\alpha
\end{equation*}
and we are now "at" the point $\rho+1$ in $\mathbb{C}$, meaning that the new factor of $\left(\frac{2\pi in}{\rho+1}\right)$ has been attached to the point $\rho+1$. Integrating by parts again deposits another factor of $\left(\frac{2\pi in}{\rho+2}\right)$, this time attached to the point $\rho+2$ (since the remaining integral is now against $\overset{\lor}{q}_{\rho+2}(\alpha)$). If we continue in the same fashion we deposit further factors of $\left(\frac{2\pi in}{\rho+j}\right)$ at $\rho+j$ for $j=3,4, \ldots$ and our formula for $a_{n}(\rho)$ becomes a product of two limiting quantities. The first is the quotient of remainder products given by $\prod_{R,+}[2\pi in](\rho)$ over $\prod_{R,+}[\tilde{z}](\rho)$. The second is $\lim_{k\rightarrow\infty}\int_{0}^{1}\overset{\lor}{q}_{\rho+k}(\alpha) e^{-2\pi in\alpha}\,d\alpha$.

Now $2\pi in$ is a constant factor and we derived in [2] that therefore, in a generalised geometric C\'{e}saro sense, we have $\prod_{R,+}[2\pi in](\rho)=(2\pi in)^{-\rho-\frac{1}{2}}$. And likewise we saw from our generalised geometric C\'{e}saro definition of the Gamma function in [2] that we have $\prod_{R,+}[\tilde{z}](\rho)=\frac{(2\pi)^{\frac{1}{2}}}{\Gamma(\rho+1)}$. Hence we have
\begin{equation}
a_{n}(\rho)=\begin{cases}
\frac{\Gamma(\rho+1)}{(2\pi)^{\rho+1}} \cdot \frac{\textrm{e}^{-i\frac{\pi}{2}(\rho+\frac{1}{2})}}{n^{\rho+\frac{1}{2}}} \cdot \lim_{k\rightarrow\infty}\int_{0}^{1}\overset{\lor}{q}_{\rho+k}(\alpha) e^{-2\pi in\alpha}\,d\alpha & \quad \textrm{if} \; n\in\mathbb{Z}\setminus\{0\}\\
0 & \quad \textrm{if} \; n=0 \quad .\end{cases}
\label{q_cech_FS_coefficients_1_alt}\end{equation}
Overall, comparing with equation \ref{q_cech_FS_coefficients_1} we see that most of the components of the formula for $a_{n}(\rho)$ have simply dropped out naturally from our geometric C\'{e}saro approach. All that remains is to prove, for the limiting integral, that
\begin{equation}
\lim_{k\rightarrow\infty}\int_{0}^{1}\overset{\lor}{q}_{\rho+k}(\alpha) \textrm{e}^{-2\pi in\alpha}\,d\alpha = \textrm{e}^{2\pi i\rho} \cdot \textrm{e}^{\frac{3\pi i}{4}} \cdot \frac{1}{\sqrt{n}} \quad .
\end{equation}
We have not succeeded in this final step to this point, but we believe deriving this limiting integral would be interesting in its own right. And we likewise believe that the above working shows that this entire geometric C\'{e}saro approach is promising, with potential application not just for this calculation but much more broadly. 

\section{Acknowledgements}

We thank Professor R.W.Emerson (pers. comm.) for many insightful comments and Professor T. Abby for his help in preparing this paper.

\end{document}